\documentclass[11pt]{article}

\usepackage{amsmath}
\usepackage{t-angles}

\title{ An Example of Double Cross Coproducts with Non-trivial
Left Coaction and Right Coaction in Strictly Braided Tensor
Categories
 \thanks {This work was supported by the  National Natural Science Foundation  (No. 19971074)}}
\author{
Shouchuan Zhang \ \ \ Bizhong Yang\\ Department  of Mathematics,
Hunan University\\ Changsha  410082, \
 P.R.China. \ \
E-mail:z9491@yahoo.com.cn\\
Beishang Ren\\
Department  of Mathematics, Guangxi Normal College\\
Nanning  530001, \
 P.R.China. \ \
}
\date{}
\begin{document}
\newtheorem{Theorem}{\quad Theorem}[section]
\newtheorem{Proposition}[Theorem]{\quad Proposition}
\newtheorem{Definition}[Theorem]{\quad Definition}
\newtheorem{Corollary}[Theorem]{\quad Corollary}
\newtheorem{Lemma}[Theorem]{\quad Lemma}
\newtheorem{Example}[Theorem]{\quad Example}
\maketitle \addtocounter{section}{-1}

 \begin {abstract}  An example of double cross coproducts with both non-trivial
left coaction and non-trivial right coaction in strictly braided
tensor categories is given.
\\
\noindent 2000 Mathematics subject Classification: 16w30.\\
Keywords: Hopf algebra, braided tensor category, double cross
coproduct.
 \end {abstract}

\section { Introduction and Preliminaries }  The double cross
coproducts in braided tensor categories have been studied by
Y.Bespalov, B.Drabant and author in \cite {BD} \cite {ZC}.
However, hitherto any examples of double cross coproducts with
both non-trivial left coaction and non-trivial right coaction in
strictly braided tensor categories (i.e. the braiding is not
symmetric ) have not been found. Therefore Professor S.Majid asked
if there is such example.

In this paper we first give the cofactorisation theorem of Hopf
algebras in braided tensor categories.  Using the cofactorisation
theorem and Sweedler four dimensional Hopf algebra, we construct
such example.

We denote the multiplication, comutiplication, evaluation $d$,
coevaluation $b$,
  braiding and inverse braiding by \[
\begin{tangle}
\cu \step[1],\step[2] \cd \step[2],  \ev \ \ , \coev\ \ ,
\step[2]\x \step[2] \hbox { and } \step[2] \xx \step[2]\ \ ,
\end{tangle}
\]
respectively. For convenience, we denote the inverse of  morphism
$f$
 by $\overline{f}$ if $f$ has an inverse.
 \newpage
 Since every braided tensor category is always
equivalent to a strict braided tensor category by \cite [Theorem
0.1] {ZC},
 we can view every braided tensor category as
 a strict braided tensor category and use braiding
diagrams freely.

\section {The cofactorisation theorem of bialgebras in braided tensor
categories }

Throughout this section, we work in braided tensor category
$({\cal C}, C)$ and assume that all  Hopf algebras and bialgebras
are living in $({\cal C}, C)$ unless otherwise stated.  We give
the cofactorisation theorem of bialgebras in braided tensor
categories in this section.

We first recall the double bicrossproducts in \cite {ZC}.  Let $H$
and $A$ be
 two bialgebras in  braided tensor categories and \ \ \ \ \
 \begin {eqnarray*}
\alpha : H \otimes A \rightarrow &A& , \hbox { \ \ \ \ }
\beta : H \otimes A \rightarrow H,    \\
\phi :  A \rightarrow H \otimes  &A&   , \hbox { \ \ \ \ } \psi :
H  \rightarrow H \otimes A
\end {eqnarray*}
morphisms in ${\cal C}$.
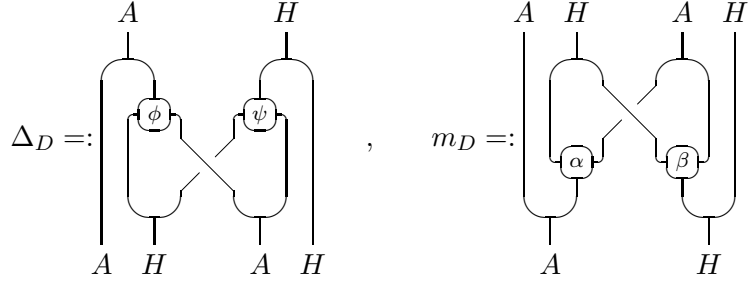
\[ \Delta _D =:
\begin{tangle}
\step\object{A}\Step\Step\Step\object{H}\\
 \cd \step[4] \cd \\
 \id \step [1]\td \phi \step [2]\td \psi \step [1] \id\\
 \id \step \id \step [2] \x \step [2]\id \step \id \\
 \id \step \cu \step [2] \cu \step [1]\id  \\
\object{A}\step\step\object{H}\Step\Step\object{A}\step\step\object{H}
 \end{tangle}
 \Step,\Step m_D =:
 \begin{tangle}
\object{A}\step\step\object{H}\Step\Step\object{A}\step\step\object{H}\\
 \id \step \cd \step [2] \cd \step [1]\id  \\
\id \step \id \step [2] \x \step [2]\id \step \id \\
 \id \step [1]\tu \alpha \step [2]\tu \beta \step [1] \id\\
\cu\Step\Step\cu\\
 \step\object{A}\Step\Step\Step\object{H}
 \end{tangle}
 \]
 \noindent and $\epsilon _D = \epsilon _A
\otimes \epsilon _H$ , $ \eta _D = \eta _A \otimes \eta _H. $ We
denote   $(A \otimes H, m_D, \eta _D, \Delta _D, \epsilon _D) $
  by                      $$ A {}^{\phi} _{\alpha}  \bowtie
^{\psi}_{\beta} H ,$$ which is
 called the
double bicrossproduct of $A$ and $H$.

When $\phi $  and $\psi$ are trivial,  we   denote $ A{}_\alpha
^\phi \bowtie _\beta ^\psi H $  by
 $A{} _\alpha \bowtie  _\beta H$.  When
 $\alpha $ and $\beta $ are trivial, we denote
  $ A{}_\alpha ^\phi \bowtie _\beta ^\psi H $  by
$ A{} ^\phi \bowtie  ^\psi H  $.  We call $ A{}_\alpha  \bowtie
_\beta H $  a double cross product and denote it by  $A  {\bowtie
}  H$ in short.  We call
  $ A {}^\phi \bowtie  ^\psi H $  a double cross coproduct.

 \begin {Theorem} \label {1.1}  (Factorisation theorem) (See
\cite [Theorem 7.2.3]{Ma95b}) Let $X$ , $A$ and  $H$  be
bialgebras or  Hopf  algebras. Assume that $j_A$ and $j_H$  are
bialgebra or Hopf algebra morphisms from $A$ to $X$ and $H$  to
$X$ respectively. If $   \xi =:  m_X (j_A \otimes j_H)$ is an
isomorphism  from $ A \otimes H$ onto $X$ as objects in ${\cal
C}$, then there exist morphisms
 $$ \alpha  : H\otimes  A \rightarrow A \hbox { \ \ \ and \ \ \ }
 \beta  :  H \otimes A \rightarrow H $$
 such that $ A {}_ \alpha \bowtie _\beta H $  becomes  a bialgebra or
 Hopf algebra
  and  $\xi $   is a bialgebra or Hopf algebra isomorphism from $ A {}_\alpha \bowtie _\beta H$  onto
  $X.$
\end {Theorem}
\textbf{Proof}. Set
\[  \zeta  =:
\begin{tangle}
\step\object{H}\Step\object{A}\\
 \step \O {j_H} \step[2] \O {j_A} \\
\step \cu \\
\step \td {\bar \xi}\\
\step \object{A}\step\step\object{H}
 \end{tangle}
 \Step,\Step \alpha =:
 \begin{tangle}
\object{H}\step\step\object{A}\\
 \ox \zeta \\
 \id \step[2] \QQ \epsilon \\
 \object{A}
 \end{tangle}\step \hbox { \ and } \step
 \beta =:
 \begin{tangle}
\step \object{H}\step\step\object{A}\\
 \step \ox \zeta \\
\step \QQ \epsilon \step[2] \id  \\
 \step [3]\object{H} \ \ \
 \end{tangle}
 \ \ .\]
We see
\[
\begin{tangle}
\object{H}\Step\object{H}\step\step\object{A}\\
\id\step\step\ox \zeta \\
        \ox \zeta\step\step\id \\
              \O {j_{A}}\step\step\cu \\
       \id\step\step\step\O {j_{H}}\\
       \id\step\step\dd\\
       \cu\\
       \step\object{X}\end{tangle}
\;=\enspace
\begin{tangle}
\object{H}\Step\object{H}\step\step\object{A}\\
\id\step\step\ox \zeta \\
        \ox \zeta \step\step\id \\
        \O{j_{A}}\Step\O{j_{H}}\Step\O{j_{H}}\\
       \id\Step\cu \\
       \id\Step\dd \\
       \cu \\
       \step\object{X}\end{tangle}
\;=\enspace
\begin{tangle}
\object{H}\Step\object{H}\step\step\object{A}\\
\id\step\step\ox \zeta \\
        \ox\zeta \step\step\id \\
        \O{j_{A}}\Step\O{j_{H}}\Step\O{j_{H}}\\
       \cu\Step\id \\
       \step\d\Step\id \\
       \Step\cu \\
       \Step\step\object{X}\end{tangle}
       \;=\enspace
\begin{tangle}
\object{H}\Step\object{H}\step\step\object{A}\\
\id\step\step\ox \zeta \\
  \O{j_{H}}\Step\O{j_{A}}\Step\O{j_{H}}\\
       \cu\Step\id \\
       \step\d\Step\id \\
       \Step\cu\\
       \Step\step\object{X}\end{tangle}
  \;=\enspace
\begin{tangle}
\object{H}\Step\object{H}\step\step\object{A}\\
\id\step\step\id\Step\id\\
  \O{j_{H}}\Step\O{j_{H}}\Step\O{j_{A}}\\
       \cu\Step\id \\
       \step\d\Step\id \\
       \Step\cu \\
       \Step\step\object{X}\end{tangle}
\ =\
\begin{tangle}
\object{H}\Step\object{H}\step\step\object{A}\\
\cu\Step\id\\
\step\id\Step\dd\\
\step \ox \zeta \\
\step\O{j_{A}}\Step\O{j_{H}}\\
\step\cu \\
\step\step\object{X} \ \ \end{tangle} \ . \] Thus
\[
\begin{tangle}
\object{H}\Step\object{H}\step\step\object{A}\\
\id\step\step\ox \zeta \\
\ox \zeta \Step\id\\
\id\Step\cu\\
\object{A}\Step\step\object{H}
\end{tangle}
\;=\enspace
\begin{tangle}
\object{H}\Step\object{H}\step\step\object{A}\\
\cu \Step\id\\
\step\id\Step\dd\\
\step \ox \zeta \\
\step\object{A}\Step\object{H}\end{tangle}\ \ . \Step\Step \cdots
\cdots  (1)
\]
Similarly we have
\[
\begin{tangle}
\object{H}\Step\object{A}\step\step\object{A}\\
\ox \zeta \Step\id\\
\id\Step\ox \zeta \\
\cu \Step\id\\
\step\object{A}\Step\object{H}
\end{tangle}
\;=\enspace
\begin{tangle}
\object{H}\Step\object{A}\step\step\object{A}\\
\id\Step\cu\\
\d\Step\id\\
\step \ox \zeta \\
\step\object{A}\Step\object{H}\end{tangle} \ \ . \Step\Step \cdots
\cdots  (2)
\]
We also have
\begin {eqnarray*} \zeta (\eta \otimes id ) = id \otimes \eta
\hbox { \ and \ } \zeta (id \otimes \eta  ) = \eta\otimes id \  .
\ \ \ \ \ \ \ \ \cdots \cdots (3)
\end {eqnarray*}

 It is clear that $\zeta$ is a coalgebra morphism from
$H\otimes A$ to $A\otimes H$, since $j_{A}$,  $ j_{H}$ and $m_{X}$
all are coalgebra homorphisms. Thus we have
\[
\begin{tangle}
\step\object{H}\Step\Step\object{A}\\
\cd \Step\cd \\
\id\Step\x\Step\id\\
\ox \zeta \Step\ox \zeta \\
\object{A}\Step\object{H}\Step\object{A}\Step\object{H}
\end{tangle}
\;=\enspace
\begin{tangle}
\Step\object{H}\Step\object{A}\\
\Step\ox \zeta \\
\step\dd\Step\d\\
\cd\Step\cd \\
\id\Step\x\Step\id \\
\object{A}\Step\object{H}\Step\object{A}\Step\object{H}
\end{tangle}
\step \hbox { and \ } (\epsilon \otimes \epsilon )\zeta =
(\epsilon \otimes \epsilon ) . \ \ \ \ \ \ \ \ \ \ \cdots \cdots
(4)
\]
We now show that $(A,\alpha)$ is an $H$-module coalgebra:\\
\[
\begin{tangle}
\object{H}\Step\object{H}\step\object{A}\\
\cu \step\id\\
\step\tu \alpha\\
\Step\object{A}
\end{tangle} \ \
= \begin{tangle}
\object{H}\Step\object{H}\step\step\object{A}\\
\cu \Step\id\\
\step\id\Step\dd\\
\step \ox \zeta \\
\step\object{A}\Step\object{\HH\obox 1\epsilon}
\end{tangle}
 \Step \stackrel { \hbox{ by } (1)}{= } \Step
\begin{tangle}
\object{H}\Step\object{H}\step\step\object{A}\\
\id\step\step\ox \zeta \\
\ox \zeta \Step\id\\
\id\Step\cu\\
\object{A}\Step\step\object{\HH\obox 1\epsilon}
\end{tangle}
 \;=\enspace
\begin{tangle}
\object{H}\step\object{H}\Step\object{A}\\
\id\step\tu \alpha\\
\tu \alpha\\
\step\object{A}
\end{tangle}
\]
and \ \  $\alpha (\eta \otimes id_A ) = (id_A \otimes \epsilon )
\zeta (\eta \otimes id_A ) \stackrel {\hbox {by} (3)}{=} id_A .$

We see that  $\epsilon \circ \alpha  = (\epsilon \otimes \epsilon
)\zeta \stackrel {\hbox {by} (4)}{=}\epsilon \otimes \epsilon $
 \ \  and \[
\begin{tangle}
\step\object{H}\Step\Step\object{A}\\
\cd \Step\cd\\
\id\Step\x\Step\id\\
\cu\Step\cu\\
\step\object{A}\Step\Step\object{A}
 \end{tangle}
 \;=\enspace
\begin{tangle}
\step\object{H}\Step\Step\object{A}\\
\cd \Step\cd \\
\id\Step\x\Step\id\\
\ox \zeta \Step\ox \zeta \\
\object{A}\Step \object{\HH\obox
1\epsilon}\Step\object{A}\Step\object{\HH\obox 1\epsilon}\Step
\end{tangle}
\ \ \stackrel {\hbox {by} (4)}{= } \ \
\begin{tangle}
\Step\object{H}\Step\object{A}\\
\Step\ox \zeta \\
\step\dd\Step\d\\
\cd \Step\cd \\
\id\Step\ox \zeta \Step\id\\
\object{A}\Step\object{\HH\obox
1\epsilon}\Step\object{A}\Step\object{\HH\obox 1\epsilon}
\end{tangle}
\;=\enspace
\begin{tangle}
\object{H}\Step\object{A}\\
\cu\\
\cd\\
\object{A}\Step\object{A}
\end{tangle} \ \ .
\]Thus $(A,\alpha)$ is an $H$-module
coalgebra. Similarly, we can show that $(H,\beta)$ is an
$A$-module coalgebra.

  Now we show that conditions $(M1)$--$(M4)$ in [12,p37] hold. By (3), we easily
know that$(M1)$ holds. Next we show that$(M2)$ holds.
\[
\begin{tangle}
\step\object{H}\Step\Step\object{A}\Step\step\object{A}\\
\cd \Step\cd \Step\id\\
\id\Step\x\Step\id\Step\id\\
\tu \alpha \Step\tu \beta \Step\id\\
\step\d\Step \step \d\Step\id\\
\step\step\d\step\step\step\tu \alpha\\\
\step\step\step\Cu\\
\Step\Step\step\object{A}
\end{tangle}
\;=\enspace
\begin{tangle}
\step\object{H}\Step\Step\object{A}\Step\step\object{A}\\
\cd \Step\cd \Step\id\\
\id\Step\x\Step\id\Step\id\\
\ox \zeta \Step\ox \zeta \Step\id\\
\id\Step\QQ \epsilon \Step\QQ \epsilon \Step\ox \zeta  \\
\d\Step\Step\dd \Step \QQ \epsilon \\
\step\Cu\\
\Step\step\object{A}\Step
\end{tangle}
\ \ \stackrel {\hbox {by } (4)} {= } \
\begin{tangle}
\Step\object{H}\Step\object{A}\Step\Step\object{A}\\
\Step\ox \zeta \Step\Step\id\\
\step\dd\Step\d\Step\step\id\\
\cd \Step\cd \Step\id\\
\id\Step\ox \zeta \Step\id\Step\id\\
\id\Step\QQ \epsilon\Step\QQ \epsilon\Step\ox \zeta \\
\d\Step\Step\dd\Step\QQ \epsilon\\
\step\Cu\\
\Step\step\object{A}
\end{tangle}
\]
\[
\;=\enspace
\begin{tangle}
\object{H}\Step\object{A}\step\step\object{A}\\
\ox \zeta \Step\id\\
\id\Step\ox \zeta \\
\cu \Step\id\\
\step\object{A}\Step\step\object{\HH\obox 1\epsilon}
\end{tangle}
\;=\enspace
\begin{tangle}
\object{H}\Step\object{A}\step\step\object{A}\\
\nw1 \step\cu\\
\step \ox \zeta \\
\step \id \step [2] \QQ \epsilon \\
 \step\object{A}\end{tangle} \step \stackrel
{\hbox {by} (2)}{=} \step
\begin{tangle}
\object{H}\Step\object{A}\step\step\object{A}\\
\id\Step\cu\\
\d\Step\id\\
\step\tu \alpha\\
\Step\object{A}
\end{tangle} \ .
\]
 Thus $(M2)$ holds. Similarly, we can get the proofs of
$(M3)$ and$(M4)$. Consequently,   $ A {}_\alpha \bowtie _\beta H$
is a bialgebra or Hopf algebra by [12,  Corollary 1.8]. It
suffices to show that $\zeta$ is a bialgebra morphism from$ A {}
_\alpha \bowtie _\beta H$ to $X$. Let $ D=A {}_\alpha \bowtie
_\beta H$. Since
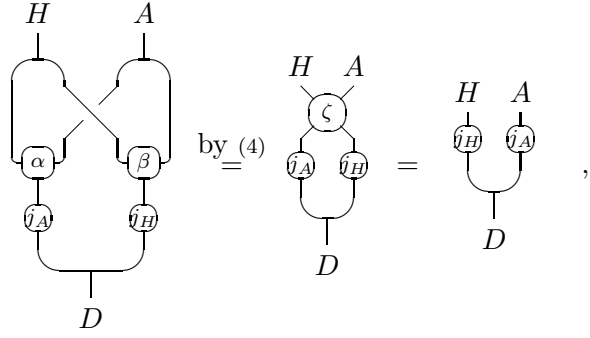
\[
\begin{tangle}
\step\object{H}\Step\Step\object{A}\\
\cd\Step\cd\\
\id\Step\x\Step\id\\
\tu \alpha \Step\tu \beta\\
\morph {j_{A}} \Step\morph {j_{H}}\\
\step\Cu\\
 \step\Step\object{D}
\end{tangle}
\ \ \stackrel {\hbox {by } (4)} {= }
\begin{tangle}
\step\object{H}\Step\object{A}\\
\step\ox \zeta \\
\morph {j_{A}} \morph {j_{H}}\\
\step\cu\\
\Step\object{D}
\end{tangle}
\;=\enspace
\begin{tangle}
\step\object{H}\Step\object{A}\\
\morph {j_{H}} \morph {j_{A}}\\
\step\cu\\
\Step\object{D}
\end{tangle}
\step\;,  \enspace
\]
 we have that $\xi$ is a bialgebra morphism from $ A{} _\alpha \bowtie _\beta H$ to $X$ by \cite
 [Lemma 2.5] {ZC}. \begin{picture}(5,5)
\put(0,0){\line(0,1){5}}\put(5,5){\line(0,-1){5}}
\put(0,0){\line(1,0){5}}\put(5,5){\line(-1,0){5}}
\end{picture}\\

\begin {Theorem} \label {1.2}
(Co-factorisation theorem) Let $X$ ,   $A$ and  $H$  be bialgebras
or  Hopf  algebras. Assume that $p_A$ and $p_H$  are bialgebra or
Hopf algebra morphisms from $X$ to $A$ and $X$  to $H$,
respectively. If  $ \xi =   (p_A \otimes p_H)\Delta _X$ is an
isomorphism  from $X$ onto $ A \otimes H$ as objects in ${\cal
C}$,   then there exist morphisms:
 $$ \phi  :  A \rightarrow H\otimes A \hbox { \ \ \ and \ \ \ }
 \psi  :  H \rightarrow H \otimes A$$
 such that $ A {}^\phi \bowtie ^\psi H $  becomes  a bialgebra or Hopf
 algebra   and  $\xi $   is a bialgebra or Hopf algebra isomorphism from
  $X$ to $ A {}^\phi \bowtie ^\psi H$.
  \end {Theorem}

  {\bf Proof.}  Set \\
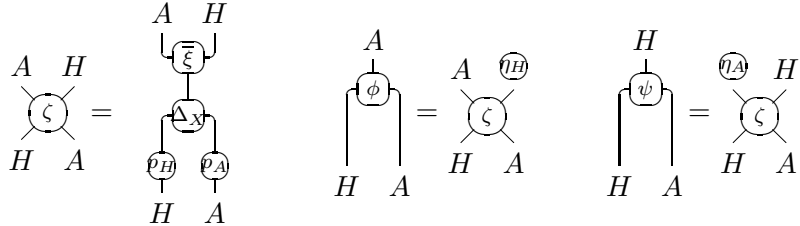
\[
\begin{tangle}
\object{A}\Step\object{H}\\
\ox \zeta \\
\object{H}\Step\object{A}
\end{tangle}
\;=\enspace
\begin{tangle}
\step\object{A}\Step\object{H}\\

\step\tu {\overline{\xi}}\\
\step\td {\Delta_{X}} \\
\morph {p_{H}} \morph {p_{A}}\\
\step\object{H}\Step\object{A}
\end{tangle}
\Step\Step
\begin{tangle}
\step\object{A}\\
\td \phi \\
\id\step\step\id\\
\object{H}\Step\object{A}
\end{tangle}
\;=\enspace
\begin{tangle}
\object{A}\Step\object{\HH\obox 1{\eta_{H}}}\\
\ox \zeta \\
\object{H}\Step\object{A}\\
\end{tangle}
\Step\Step
\begin{tangle}
\step\object{H}\\
\td \psi \\
\id\step\step\id\\
\object{H}\Step\object{A}
\end{tangle}
\;=\enspace
\begin{tangle}
\object{\HH\obox 1{\eta_{A}}}\Step\object{H}\\
\ox \zeta \\
\object{H}\Step\object{A}\\
\end{tangle}
\] We can complete the proof by turning upside down the diagrams in
 the proof of the preceding theorem. \begin{picture}(5,5)
\put(0,0){\line(0,1){5}}\put(5,5){\line(0,-1){5}}
\put(0,0){\line(1,0){5}}\put(5,5){\line(-1,0){5}}
\end{picture}\\

 From now on,   we always consider Hopf algebras over field $k$ and
 the diagram
 
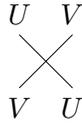
\[
 \begin{tangle}
 \object{U}\Step\object{V}\\
 \XX \\
 \object{V}\Step\object{U}
 \end{tangle}
 \]
 always denotes the ordinary twisted map:
 $x\otimes y \longrightarrow y\otimes x$.
 Our diagrams only denote homomorphisms between vector spaces,   so
two diagrams can have the additive operation.

Let $H$  be an ordinary bialgebra and $(H_1,   R)$   an ordinary
quasitriangular Hopf algebra over field $k$. Let  $f $ be
 a  bialgebra homomorphism from $H_1 $  to $H$.
 Then there exists  a   bialgebra
 $B$,   written as  $B(H_1 ,   f,   H)$,    living in $({}_{H_1} {\cal M},   C^R).$
  Here
  $B(H_1,   f,   H) = H$ as algebra,
  its counit is $\epsilon _H$,    and its comultiplication and antipode
are\\
\[
\begin{tangle}
\step \object{B}\\
\td {\Delta_{B}}\\
\id \Step \id\\
\object{H}\Step\object{H}\\
\end{tangle}
\ = \
\begin{tangle}
\step [1]\Step\object{B}\\
\step [2] \td {\Delta_{H}}\\

\step \ne1 \step [2] \nw1\\
 \dd \step \ro R\step \d\\
\id \step\morph f \morph {fS}\step \id\\
\id \step\step \XX  \step\step \id \\
\cu  \Step   \tu {ad}\\
\step\object{H} \Step\step\step \object{H}
\end{tangle} \ \ , \ \
\Step\Step\Step
\begin{tangle}
\object{B}\\
\O {S_{B}}\\
 \object{B}
\end{tangle}
\step\;=\enspace
\begin{tangle}
\Step\Step\Step\object{B}\\
\step  \ro R  \Step\step \id\\
 \morph f  \morph f \step \step \id\\
\step \nw3 \step\id\Step \step\id\\
\step\Step\id\step\tu {ad}\\
\step\Step\id\step\morph S\\
\step\Step \cu \\
\Step\Step\object{B}
 \end{tangle}
\]
 (see \cite [Theorem 4.2]{Ma95a}), respectively.
  In particular,   when $H=H_1$ and $f= id _H$,    $B(H_1,   f,   H)$ is  a
  braided group,   called
  the braided group analogue of $H$ and written as $\underline H$.

$R$ is called a weak $R$-matix of  $A\otimes H$ if $R$ is
invertible under convolution with
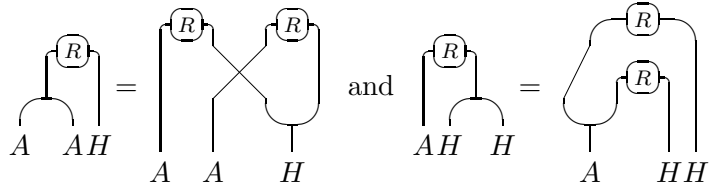
\[
\begin{tangle}
\step \ro R  \\
 \cd \step \id\\
\object{A}\Step\object{A}\step\object{H}\\
\end{tangle}
\;= \enspace
\begin{tangle}
\ro R  \Step \ro R  \\
\id \Step \XX \Step\id\\
\id \Step \id \Step \cu\\
\object{A}\Step\object{A}\Step\step\object{H}
\end{tangle} \hbox { \ \ and \ \ }
\begin{tangle}
\ro R  \\
\id \step \cd\\
\object{A}\step\object{H}\Step\object{H}
\end{tangle}
\;= \enspace
\begin{tangle}
\step \Ro R\\
\dd \step \ro R \step  \id\\
\cu \Step \id \step \id\\
\step\object{A}\Step\step\object{H}\step\object{H}
\end{tangle}
\]

Let $(A,  P)$ and $(H,  Q)$ be ordinary finite-dimensional
quasitriangular Hopf algebras over field $k$. Let $R$ be a weak
$R$-matrix of $A\otimes H$. For any $U,  V \in CW(A\otimes
H)=:$$\{U \in A\otimes H \ | \ U$  is a weak $R$-matrix and in the
center of $ A\otimes H\}$,
 $$R_{D}=:\sum R'P'U' \otimes Q'(R^{-1})''V'' \otimes P'' (R^{-1})' V'
 \otimes R''Q''U''$$
is a quasitriangular structure of $D$ and every quasitriangular
structure of $D$ is of this form (\cite [Theorem 2.9 ] {Ch}),
where $R = \sum R'\otimes R''$,   etc.

\begin {Lemma} \label {1.3}
Under the above discussion,   then

(i) $\pi _A : D \rightarrow A$ and $\pi _H : D \rightarrow H$ are
bialgebra  or Hopf algebra homomorphisms,   respectively. Here
$\pi _A $ and $\pi _H$ are trivial action,   that is,   $\pi _A (h
\otimes a ) = \epsilon (h)a$ for any $a\in A,   h \in H.$

(ii) $B(D,   \pi _A,   A) = \underline A$ and $B(D,   \pi _H,   H)
= \underline H.$

(iii) $\pi _{\underline A} : \underline D \rightarrow \underline
A$ and $\pi _{\underline H }: \underline D \rightarrow \underline
H$ are bialgebra  or Hopf algebra homomorphisms,   respectively.

  \end {Lemma}
  {\bf Proof.} (i) It is clear.

  (ii) It is enough to show $\Delta _B = \Delta _{\underline
  {A}}$ since $B= \underline A$ as algebras,   where $B =: B(D,   \pi_ A,
  A).$  See\\
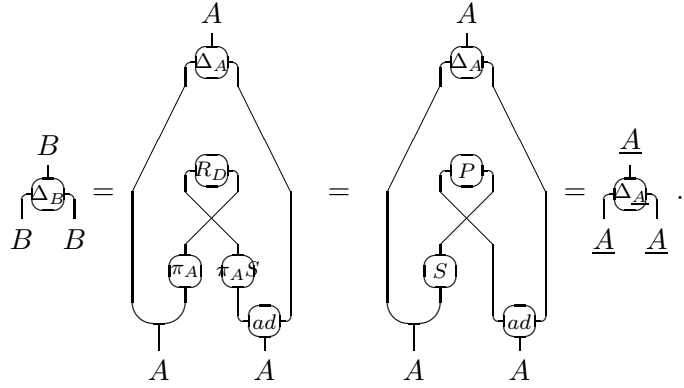
\[
\begin{tangle}
\step \object{B}\\
\td {\Delta_{B}} \\
\object{B}\Step\object{B}\\
\end{tangle}
\ = \
\begin{tangle}
\step\Step\object{A}\\

\step [2] \td {\Delta_{A}} \\

\step \ne1 \step [2] \nw1\\
 \dd \step \ro {R_{D}}\step \d\\
\id\step\step\XX  \step\step\id \\
\id \Step \O {\pi_{A}} \Step\O {\pi_{A}S}\Step\id\\
\cu \step  \step \tu {ad}\\
\step\object{A}\Step\Step\object{A}
 \end{tangle}
\ \ \ = \ \ \
\begin{tangle}
\step\Step\object{A}\\

\step [2] \td {\Delta_{A}} \\

\step \ne1 \step [2] \nw1\\
 \dd \step \ro {P} \step \d\\
\id\step\step\XX  \step\step\id \\
\id\Step \S\Step\id\Step\id\\
\cu \step  \step \tu {ad}\\
\step\object{A}\Step\Step\object{A}
 \end{tangle}
\ = \
 \begin{tangle}
\step \object{\underline{A}}\\
\td {\Delta_{\underline{A}}} \\
\object{\underline{A}}\Step\object{\underline{A}}\\
\end{tangle} \ \ .
\]
  Thus $\Delta _B = \Delta _{\underline A}$.
Similarly, we have $B(D,   \pi _H,   H)= \underline H.$

 (iii) See\\
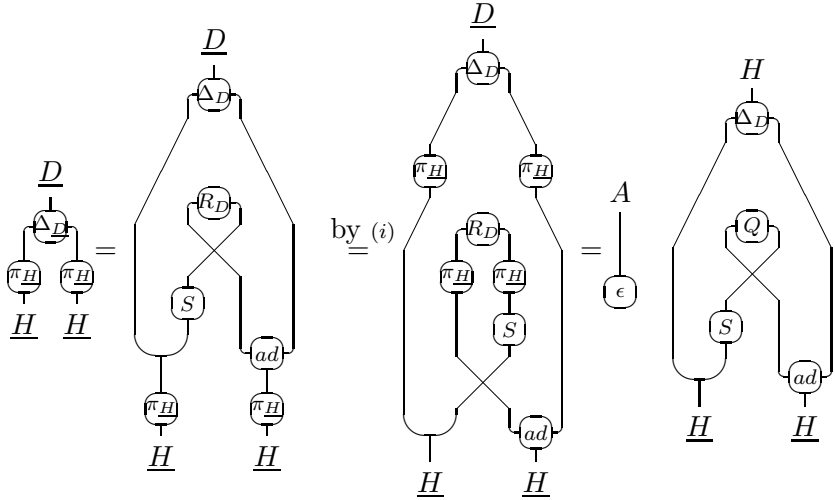
  \[
 \begin{tangle}
\step \object{\underline{D}}\\
\td {\Delta_{\underline{D}}} \\
\O {\pi_{\underline{H}}}\Step\O {\pi_{\underline{H}}}\\
\object{\underline{H}}\Step\object{\underline{H}}\\
\end{tangle}
\ = \
\begin{tangle}
\step\Step\object{\underline{D}}\\

\step [2] \td {\Delta_{D}} \\

\step \ne1 \step [2] \nw1\\

 \dd \step \ro {R_{D}}\step \d\\
\id\step\step\XX  \step\step\id \\
\id\Step \S\Step\id\Step\id\\
\cu \step  \step \tu {ad}\\
\step\O {\pi_{\underline{H}}}\Step\Step\O {\pi_{\underline{H}}}\\
\step\object{\underline{H}}\Step\Step\object{\underline{H}}
 \end{tangle}
\ \ \stackrel {\hbox { by }(i)}{ = }
\begin{tangle}
\Step \step\object{\underline{D}}\\

\step [2] \td {\Delta_{D}} \\

\step \ne1 \step [2] \nw1\\
\step\O {\pi_{\underline{H}}}\Step\step\step\O {\pi_{\underline{H}}}\\
\dd \step  \ro {R_{D}} \step \d\\
\id \Step \O {\pi_{\underline{H}}}\Step\O {\pi_{\underline{H}}}\Step\id\\
\id\Step\id\Step\S\Step\id\\
\id\step\step\XX  \step\step\id \\
\cu \step  \step \tu {ad}\\
\step\object{\underline{H}}\Step\Step\object{\underline{H}}
 \end{tangle}
   \ = \
\begin{tangle}
\object{A}\\
\id\\
\QQ \epsilon
 \end{tangle}
 \Step
 \begin{tangle}
\step\Step\object{H}\\

\step [2]\td {\Delta_{D}} \\

\step \ne1 \step [2] \nw1\\
 \dd \step \ro {Q}\step \d\\
\id\step\step\XX \step\step\id \\
\id\Step \S\Step\id\Step\id\\
\cu \step  \step \tu {ad}\\
\step\object{\underline{H}}\Step\Step\object{\underline{H}}
 \end{tangle}
\]
 and $\epsilon \circ \pi _{\underline H}= \epsilon.$ Thus $\pi
_{\underline H}
 $ is a coalgebra homomorphism.

 Since the multiplications in $\underline D$ and $\underline H$
 are the same as in $D$ and $H$, respectively,   we have that $\pi _{\underline H}
 $ is an algebra homomorphism by  (i). Similarly,   we can show
 that $\pi _{\underline A}$ is a bialgebra homomorphism.  \begin{picture}(5,5)
\put(0,0){\line(0,1){5}}\put(5,5){\line(0,-1){5}}
\put(0,0){\line(1,0){5}}\put(5,5){\line(-1,0){5}}
\end{picture}\\

 We now investigate the relation among
 braided group analogues of  quasitriangular Hopf algebras $A$  and $H$
and their double cross coproduct $ D= A  \bowtie ^R H.$

 \begin {Theorem}  \label {1.4}
  Under the above discussion,   let
  $\xi = (\pi _{\underline A}
  \otimes \pi _{\underline H}) \Delta _{\underline D}.$
 Then
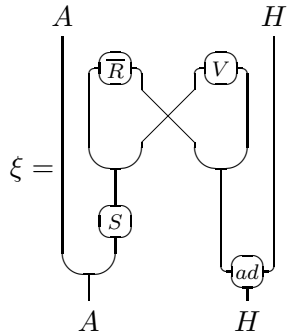
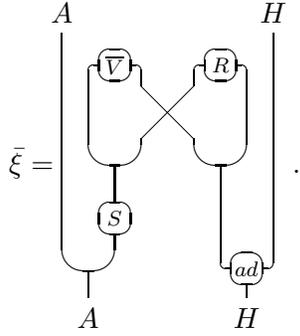
\[ \xi =
 \begin{tangle}
\object{A}\Step\Step\Step\Step \object{H}\\
\id\step\ro {\overline R} \Step \ro V
\step\id\\
\id\step\id\step\step\XX\step\step\id\step\id\\
\id\step\cu\step\step\cu\step\id\\
\id\step\step\S\step\step\step\step\id\step\step\id\\
\cu\step\step\step\step\tu {ad}\\
\step\object{A}\Step\Step\Step\object{H}
 \end{tangle}
\]
 and  $\xi $  is surjective,
where $ad$ denotes the left adjoint action of $H$.

(ii) Furthermore,   if $A$ and $H$  are finite-dimensional,   then
$\xi$  is a bijective map from $\underline D$ onto $\underline A
\otimes \underline H$. That is,
 in braided tensor category $(_{D} {\cal M},   C^{R_D}),  $
 there exist morphisms $\phi $  and $\psi $
 such that $$ \underline D \cong
 \underline A{} ^\phi \bowtie ^\psi \underline H
\hbox { \ \ \ ( as Hopf algebras )}$$ and the isomorphism is $(\pi
_{\underline A} \otimes \pi _{\underline H} ) \Delta _{\underline
D}$.

(iii) If $H$ is commutative or $V=R$,     then $\xi = id
_{\underline D} $.

\end {Theorem}
{\bf Proof.} (i)
\[
\step\step
\begin{tangle}
\object{\underline{D}}\\
\O \xi\\
\object{\underline{A}\otimes \underline{H}}
\end{tangle}
\\ \ \ \ \ = \ \ \
\begin{tangle}
\step \object{\underline{D}}\\
\td {\Delta_{\underline{D}}} \\
\O {\pi_{\underline{A}}}\Step\O {\pi_{\underline{H}}}\\
\object{\underline{A}}\Step\object{\underline{H}}\\
\end{tangle}
\ \ \ = \ \ \
\begin{tangle}
\step\Step\object{\underline{D}}\\

\step [2] \td {\Delta_{D}} \\

\step \ne1 \step [2] \nw1\\

 \dd \step\ro {R_{D}}\step \d\\
\id\Step \id\Step\S\Step\id\\
\id\step\step\XX \step\step\id \\
\cu \step  \step \tu {ad}\\
\step\O {\pi_{\underline{A}}}\Step\Step\O {\pi_{\underline{H}}}\\
\step\object{\underline{A}}\Step\Step\object{\underline{H}}
 \end{tangle}
\stackrel { \hbox {since } \pi _A, \pi _H \hbox { are algebra
homomorphisms }}{=}
\]\[
\begin{tangle}
\Step \step\object{\underline{D}}\\

\step [2] \td {\Delta_{D}} \\

\step \ne1 \step [2] \nw1\\

\step\O {\pi_{\underline{A}}}\Step\step\step\O {\pi_{\underline{H}}}\\
\dd  \step [1] \ro {R_{D}} \step \d\\
\id \Step \O {\pi_{\underline{H}}}\Step\O {\pi_{\underline{A}}}\Step\id\\
\id\Step \id\Step\S\Step\id\\
\id\step\step\XX \step\step\id \\
\cu \step  \step \tu {ad}\\
\step\object{\underline{A}}\Step\Step\object{\underline{H}}
 \end{tangle}
 \ \ \ = \ \ \
 \begin{tangle}
\Step\Step \object{\underline{D}}\\

\step\step [2] \td {\Delta_{D}} \\

\step\step \ne1 \step [2] \nw1\\
\step\step\O {\pi_{\underline{A}}}\Step\step\step\O {\pi_{\underline{H}}}\\
\step\ne2 \Step\Step\nw2\\
 \id\step\ro {\overline R}
\Step
\ro {V} \step\id\\
\id\step\id\step\step\XX\step\step\id\step\id\\
\id\step\cu\step\step\cu\step\id\\
\id\step\step\S\step\step\step\step\id\step\step\id\\
\cu\step\step\step\step\tu {ad}\\
\step\object{\underline{A}}\Step\Step\Step\object{\underline{H}}
 \end{tangle}
\ \ \ = \ \ \
 \begin{tangle}
\object{A}\Step\Step\Step\Step \object{H}\\
\id\step\ro {\overline{R}} \Step \ro {V}
\step\id\\
\id\step\id\step\step\XX\step\step\id\step\id\\
\id\step\cu\step\step\cu\step\id\\
\id\step\step\S\step\step\step\step\id\step\step\id\\
\cu\step\step\step\step\tu {ad}\\
\step\object{\underline{A}}\Step\Step\Step\object{\underline{H}}
 \end{tangle} \ \ .
\]

(ii) By the proof of  (i),   $\xi$ is bijective and\\
 \[ \bar \xi =
 \begin{tangle}
\object{A}\Step\Step\Step\Step \object{H}\\
\id\step\ro {\overline{V}} \Step \ro R
\step\id\\
\id\step\id\step\step\XX\step\step\id\step\id\\
\id\step\cu\step\step\cu\step\id\\
\id\step\step\S\step\step\step\step\id\step\step\id\\
\cu\step\step\step\step\tu {ad}\\
\step\object{A}\Step\Step\Step\object{H}
 \end{tangle} \ \ .
\]
Applying the cofactorization theorem \ref {1.2},   we complete the
proof of (ii).

(iii) follows from (i). \begin{picture}(5,5)
\put(0,0){\line(0,1){5}}\put(5,5){\line(0,-1){5}}
\put(0,0){\line(1,0){5}}\put(5,5){\line(-1,0){5}}
\end{picture}\\
{\bf Remark.} Under the assumption of Theorem \ref {1.4},   if we
set 
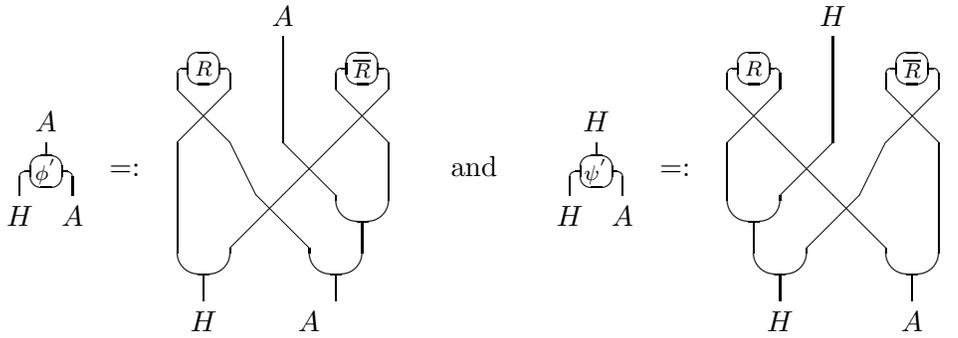
\[
 \begin{tangle}
\step\object{A}\\
\td { \phi^{'}}\\
\object{H}\Step\object{A}
 \end{tangle}
 \ \ \ =: \ \ \
  \begin{tangle}
  \Step\Step\object{A}\\
  \ro R\Step\id\Step \ro {\overline{R}}\\
  \XX\Step\id\Step\XX\\
  \id\Step\d\step\XX\Step\id\\
  \id\Step\sw2\nw2\step\step\cu\\
 \cu\Step\step\cu\\
\step\object{H}\Step\Step\object{A}
   \end{tangle}
     \Step \hbox {\ \hbox {and} \ }\Step
    \begin{tangle}
\step\object{H}\\
\td {\psi^{'}}\\
\object{H}\Step\object{A}
 \end{tangle}
 \ \ \ =: \ \ \
  \begin{tangle}
  \Step\Step\object{H}\\
  \ro R\Step\id\Step\ro {\overline{R}}\\
  \XX\Step\id\Step\XX\\
  \id\Step\XX\step\dd\Step\id\\
  \cu\step\sw2\nw2\step\step\step\id\\
\step \cu\Step\step\cu\\
\step\step\object{H}\Step\Step\step\object{A}
   \end{tangle} \ \  ,
\]
 then $A{}^{\phi '} \bowtie ^{\psi '} H= A
\bowtie ^R H$ by \cite [Lemma 1.3] {Ch}.  We now see the relation
among $\phi,   \psi,   R,   \phi'$ and $\psi '$:
\[ \phi \ \  \stackrel {\hbox {by proof of Th \ref {1.2} }}{=}\  \zeta (id \otimes \eta
)\ \  = \ \
 \begin{tangle}
\object{A}\Step\object{\HH\obox 1\eta}\\
\tu {\overline{\xi}}\\
  \td {\Delta_{\underline{D}}}\\
\O {\pi_{\underline{H}}}\Step\O {\pi_{\underline{A}}}\\
 \object{H}\Step\object{A}
\end{tangle} \ \ \ \  \stackrel {\hbox {by Th \ref {1.4} (iii)}} {=}\ \ \
 \]
 \[\begin{tangle}
\step\object{A}\Step\Step\Step\Step\object{\HH\obox 1\eta}\\
\step\id\step\ro {\overline{V}}\step\step
\ro R\step\id\\
\step\id\step\id\step\step\XX\step\step\id\step\id\\
\step\id\step\cu\step\step\cu\step\id\\
\step\id\step\step\S\step\step\step\step\id\step\step\id\\
\step\cu\Step\Step\tu {ad}\\
\step\cd\Step\Step\cd\\
\step\id\step\td {\phi^{'}}\Step
\td {\psi^{'}}\step\id\\
\step\QQ \epsilon \step\id\Step\XX\Step\id \step\QQ \epsilon\\
\step\step\cu\Step\cu\\
\Step\ne2\Step\step\step\nw2\\
\step\id\step  \ro R\Step
 \ro U \step\id\\
\step\id\step\id\Step\XX\Step\id\step\id\\
\step\id\step\cu\Step\cu\step\id\\
\step\id\Step\id\Step\Step\S\Step\id\\
\step\id\Step\nw2\step\step\ne2\Step\id\\
\step\id\Step\sw2\step\step\se2\Step\id\\
\step\cu\Step\Step\tu {ad}\\
\step\step\object{H}\Step\Step\Step\object{A}
\end{tangle}
  \ = \
 \begin{tangle}
 \Step \step\Step\object{A}\\
 \step [4]\td {\phi ^{'}}\\

\step[3] \ne1 \step [2] \nw1\\
 \Step\ne2\Step\step\step\nw2\\
\step\id\step  \ro R\Step
 \ro U\step\id\\
\step\id\step\id\Step\XX\Step\id\step\id\\
\step\id\step\cu\Step\cu\step\id\\
\step\id\Step\id\Step\Step\S\Step\id\\
\step\id\Step\nw2\step\step\ne2\Step\id\\
\step\id\Step\sw2\step\step\se2\Step\id\\
\step\cu\Step\Step\tu {ad}\\
\step\step\object{H}\Step\Step\Step\object{A}
\end{tangle}
 \ = \
 \begin{tangle}
 \Step\Step\Step\object{A}\\
 \ro R\Step\Step\id\Step\ro {\overline{R}}\\
 \XX\Step\Step\id \step[2] \XX\\

 \id\Step\nw2\Step\step[1]\XX\Step\id\\

 \id\Step\Step\XX\Step\id\step[2]\id\\

 \id\Step \Step\id\Step\nw1\step [1]\cu\\

 \id\Step \Step\id\step [3]\cu \\

 \id\Step\Step\id\Step\Step\id\Step\\
 \Cu\Step\Step \nw2\\
 \Step\id\step  \ro {R}\Step
 \ro U\step\id\\
\Step\id\step\id\Step\XX\Step\id\step\id\\
\Step\id\step\cu\Step\cu\step\id\\
\Step\id\Step\id\Step\Step\S\Step\id\\
\Step\id\Step\nw2\step\step\ne2\Step\id\\
\Step\id\Step\sw2\step\step\se2\Step\id\\
\Step\cu\Step\Step\tu {ad}\\
\Step\step\object{H}\Step\Step\Step\object{A}
 \end{tangle}
\]
 $ \hbox {and \ \ \ } \psi \ \  \stackrel { \hbox { by proof of Th \ref {1.2} }}{=}
  \ \ \  \zeta (\eta  \otimes
id  ) = $
 \[\begin{tangle}
\object{\HH\obox 1\eta}\Step\object{H}\\
\tu {\overline{\xi}}\\
 \td {\Delta _{\underline D}} \\
\O {\pi_{\underline{H}}}\Step\O {\pi_{\underline{A}}}\\
 \object{H}\Step\object{A}
\end{tangle} \ \ \  \stackrel { \hbox{ by Th \ref {1.4} (iii) } }{=} \ \ \
 \begin{tangle}
\step\object{\HH\obox 1\eta}\Step\Step\Step\Step\object{H}\\
\step\id\step\ro {\overline{V}}\step\step
\ro R\step\id\\
\step\id\step\id\step\step\XX\step\step\id\step\id\\
\step\id\step\cu\step\step\cu\step\id\\
\step\id\step\step\S\step\step\step\step\id\step\step\id\\
\step\cu\Step\Step\tu {ad}\\
\step\cd\Step\Step\cd\\
\step\id\step\td {\phi^{'}}\Step
\td {\psi^{'}}\step\id\\
\step\QQ \epsilon \step\id\Step\XX\Step\id \step\QQ \epsilon\\
\step\step\cu\Step\cu\\
\Step\ne2\Step\step\step\nw2\\
\step\id\step  \ro R\Step
 \ro U\step\id\\
\step\id\step\id\Step\XX\Step\id\step\id\\
\step\id\step\cu\Step\cu\step\id\\
\step\id\Step\id\Step\Step\S\Step\id\\
\step\id\Step\nw2\step\step\ne2\Step\id\\
\step\id\Step\sw2\step\step\se2\Step\id\\
\step\cu\Step\Step\tu {ad}\\
\step\step\object{H}\Step\Step\Step\object{A}
\end{tangle}
\]\[
 \ \ \ = \ \ \
  \begin{tangle}
  \step\Step\Step\Step\Step\object{H}\\
\step\step\ro {\overline{V}}\step\step
\ro {R}\step\id\\
\Step\id\step\step\XX\step\step\id\step\id\\
\step\step\cu\step\step\cu\step\id\\
\step\step\step\S\step\step\step\step\id\step\step\id\\
\Step\step\id\Step\Step\tu {ad}\\
\Step\step\id\Step\step\step\dd\\
\Step\td {\phi^{'}}\Step
\td {\psi^{'}}\\
\step\step \id\Step\XX\Step\id \\
\step\step\cu\step\step\cu\\
\Step\ne2\Step\step\step\nw2\\
\step\id\step  \ro R\Step
 \ro U\step\id\\
\step\id\step\id\Step\XX\Step\id\step\id\\
\step\id\step\cu\Step\cu\step\id\\
\step\id\Step\id\Step\Step\S\Step\id\\
\step\id\Step\nw2\step\step\ne2\Step\id\\
\step\id\Step\sw2\step\step\se2\Step\id\\
\step\cu\Step\Step\tu {ad}\\
\step\step\object{H}\Step\Step\Step\object{A}
 \end{tangle}
 \ \ \ = \ \ \
  \begin{tangle}
  \Step\step\Step\Step\Step\Step\object{H}\\
\Step\step\step\ro {\overline{V}}\step\step
\ro R\step\id\\
\Step\Step\id\step\step\XX\step\step\id\step\id\\
\Step\step\step\cu\step\step\cu\step\id\\
\Step\step\step\step\S\step\step\step\step\id\step\step\id\\
\Step\Step\step\id\Step\Step\tu {ad}\\
 \Step\Step\step\d \Step\Step\id\\
 \ro {\overline{R}}\Step\Step\id\Step\Step\id\Step\Step
 \ro {R}\\
 \XX\Step\Step\nw2\Step\ne2\Step\Step\XX\\
 \id\Step\nw2\Step\Step\ne2
 \nw2\Step\Step\ne2\Step\id\\
 \id\Step\Step\XX\Step\Step\XX\Step\Step\id\\
\Cu\step\step\nw2\Step\ne2\Step\Cu\\
 \Step\id\Step\Step\sw2\step\step\se2\Step\Step\id\\
 \Step\id\Step\Step\id\Step\Step\id\Step\Step\id\\
 \Step\Cu\Step\Step\Cu\\
 \Step\Step\id\step  \ro R\Step
 \ro U\step\id\\
\Step\Step\id\step\id\Step\XX\Step\id\step\id\\
\Step\Step\id\step\cu\Step\cu\step\id\\
\Step\Step\id\Step\id\Step\Step\S\Step\id\\
\Step\Step\id\Step\nw2\step\step\ne2\Step\id\\
\Step\Step\id\Step\sw2\step\step\se2\Step\id\\
\Step\Step\cu\Step\Step\tu {ad}\\
\Step\Step\step\object{H}\Step\Step\Step\object{A}
 \end{tangle}
\]

Furthermore,   $\psi = \psi '$  \ \ when $H$ is commutative or $R
= V.$ In this case,   $\underline A{} ^{\phi '}\bowtie
^{\psi'}\underline H= \underline
 {A {}^{ \phi '} \bowtie ^ {\psi '} H}
  = \underline {A \bowtie ^R H}$
  \ \  as Hopf algebras living
  in braided tensor category $( {}_D {\cal M},   C^{R_D}).$

\section {An example}

In this section,   using preceding  cofactorisation theorem,   we
give an example  of double cross coproducts with both non-trivial
left coaction and non-trivial right coaction in strictly braided
tensor categories.

Let $H^* = Hom (H,  k)$ be the dual of finite-dimensional Hopf
algebra $H$. $H^*$ can become a Hopf algebra under convolution
(\cite [Theorem 9.1.3] {Mo93}). That is,   for any $f,   g \in
H^*,   h,   h^{'}\in H$,
$$(f*g )(h)= \sum _{(h)} f (h_1)g (h_2),   \Delta _{H^*}(f) (h \otimes h^{'})= f
(hh^{'}),    S_{H^*} (f ) (h) = f (S(h)).$$ Assume $\{e_{x_i} \mid
i =1,   2,   \cdots ,   n\}$ is the dual basis of $\{x_i \mid i
=1,   2,   \cdots ,   n\}$. Define  $d_H = \left \{
\begin {array} {ll} H^* \otimes H & \rightarrow
k\\
f \otimes h & \rightarrow f (h)
\end {array} \right.$ and
$b_H = \left \{ \begin {array} {ll} k & \rightarrow H \otimes H^*
\\
 1 & \rightarrow \sum _{i =1}^n  x_i \otimes e_{x_i}
\end {array} \right.$. $d_H$ and $b_H$ are called evaluation and
coevoluation of $H$, respectively. It is clear that
\[
 \begin{tangle}
\object{H^{*}}\\
\id\Step\coev\\
\id\Step\id\Step\id\\
\ev\Step\id\Step\\
 \Step \Step\object{H^{*}}
 \end{tangle}
  \ \ \ = \ \ \
   \begin{tangle}
\object{H^{*}}\\
\id\\
\id\\
\object{H^{*}}\\
 \end{tangle}
\Step  \step\;,  \enspace\step\Step
  \begin{tangle}
 \Step \Step\object{H}\\
\coev\Step\id\\
\id\Step\id\Step\id\\
\id\Step\ev\\
\object{H}
 \end{tangle}
  \ \ \ = \ \ \
   \begin{tangle}
\object{H}\\
\id\\
\id\\
\object{H}\\
 \end{tangle}
 \Step  \step\;,  \enspace\step\Step
\]
\[
 \begin{tangle}
\object{H^{*}}\Step\object{H^{*}}\\
\cu\\
\step\object{H^{*}}
 \end{tangle}
 \ \ \ = \ \ \
   \begin{tangle}
\object{H^{*}}\Step\object{H^{*}}\\
\id\Step\id\Step\step\coev\\
\id\Step\id\Step\cd\step\id\\
\id\Step\XX\Step\id\step\id\\
\ev\Step\ev\step\id\\
\Step\Step\Step\step\object{H^{*}}
\end{tangle}
\Step and \Step
 \begin{tangle}
 \step\object{H^{*}}\\
 \cd\\
\object{H^{*}}\Step\object{H^{*}}
 \end{tangle}
 \ \ \ = \ \ \
  \begin{tangle}
\object{H^{*}}\\
\id\step\coev \Step\coev\\
\id\step\id\Step\XX \Step\id\\
\id\step\cu\Step\id\Step\id\\
\ev\Step\step\id\Step\id\\
\Step\Step\step \object{H^{*}}\Step\object{H^{*}}
 \end{tangle}
 \ \ \ .
\]
\ \ \ \  Note the multiplication and comultiplication of \ \ $H^*$
\ \  exactly are anti-multiplication and anti-comultiplication in
\cite [Proposition 2.4]{Ma95a}.
\[
 Let  \hspace{5mm} A=H^{*cop},  \Step  \tau =
 \begin{tangle}
\object{H}\Step\step\object{A}\\
\nw3\step\ne3\\
\Ev
 \end{tangle}
 \Step and \Step \textbf{$[b] =$} \Step
 \begin{tangle}
 \Q \eta \Step\Step \step\Q \eta\\
 \id\step \Coev\step\step\step\id\\
\object{A}\step\object{H}\Step\step\object{A}\step\object{H}\\
 \end{tangle} \ \ .
\]
Thus $\tau$ is a skew pairing and the Drinfeld Double $D(H)=A
\bowtie_{\tau} H$ (see [4] [5]). Furthermore,   $[b]$ is a
quasitrianglar structure of $D(H)$ by[10,  Theorem 10.3.6].
 \begin {Lemma}  \label {2.1}
 Let $H$ be a finite-dimensional Hopf algebra with $dim \ H > 1.$
 Then Drinfeld double $(D(H),   [b])$ is not triangular.
\end {Lemma}

{\bf Proof.} Assume that $x_i's$ are a basis of $H$ with $x_1 =
1_H.$ Set $x=: x_2$ and see
\[
\begin{tangle}
\Ro {\overline{[b]}}\\
 \object{x\otimes \epsilon}\Step\Step
 \object{\epsilon\otimes e_{x}}\\
 \end{tangle}
  \ \ \ = \ \ \
   \begin{tangle}
 \Q \eta \step \Coev \step [3]\Q \eta\\
 \id\step\S\step\Step \id\step[1]\id\\
\object{x}\step\object{\epsilon}\Step\step\object{\epsilon}\step\object{e_{x}}\\
 \end{tangle}
   \step\;=\enspace 0 \Step  and
\]
\[
\begin{tangle}
\ro {[b]}\\
\XX\\
 \object{x\otimes \epsilon}\step[3]
 \object{\epsilon\otimes e_{x}}\\
 \end{tangle}
  \ \ \ = \ \ \
   \begin{tangle}
 \Q \eta \Step \coev\Step\Q \eta\\
 \id\Step\XX\step\step\id\\
 \nw4\ne2\sw4\step\nw2\\
\object{x}\Step\object{\epsilon}\Step\object{\epsilon}\Step\object{e_{x}}\\
 \end{tangle}
   \step\;=\enspace 1 \ . \  \ \ \hbox {Thus } \Step
   \begin{tangle}
\Ro {\overline{[b]}}\\
 \object{D}\step\Step
 \object{D}\\
 \end{tangle}
  \step\;\neq \enspace\step
  \begin{tangle}
\ro {[b]}\\
\XX \\
 \object{D}\Step
 \object{D}\\
 \end{tangle} \ ,
\]
which implies that [$b$] is not triangular. \begin{picture}(5,5)
\put(0,0){\line(0,1){5}}\put(5,5){\line(0,-1){5}}
\put(0,0){\line(1,0){5}}\put(5,5){\line(-1,0){5}}
\end{picture}\\

 Let  us recall  Sweedler's four dimensional
Hopf algebra $H_4$. That is,   $H_4$ is a Hopf algebra  generated
$g$ and $x$ with relations
$$g^2 =1,   \ \ x^2 = 0,   \  \ xg =-gx $$
and $\Delta (x) = x \otimes 1 + g \otimes x,   \Delta (g) =
g\otimes g$,   $ \epsilon (x)=0,   \epsilon (g) =1 ,    S(x) =xg,
S(g)=g.$ Let $\{e_1,   e_g,   e_x,   e_{gx} \}$ denote the
 dual basis of $\{1,   g,   x,   gx\}$.

 \begin {Example}  \label {2.2}
 Let $H$ be Sweedler's four dimensional  Hopf algebra over field $k$ with $char \ k \not= 2.$
Let  $ D= D(H)$. Thus  $B =: D \bowtie ^{[b]}D$ is
quasitriangular, but it is not triangular by Lemma
 \ref {2.1}. Considering \cite [Theorem 2.5] {Ch},   $B$ has a  quasitriangular structure $R_B,  $
 defined in preceding Theorem \ref {1.4} with $U = V = 1 \otimes
 1$,   and $R_B$  never is triangular . Thus  $({}_B {\cal M},   C^{R_B})$ is
 a  strictly braided tensor
category by \cite [Theorem 10.4.2 (3)] {Mo93}. It follows from
Theorem \ref {1.4} (ii) that $D \bowtie ^{[b]}D \cong D{} ^\phi
\bowtie ^\psi D$ for some $\phi$ and $\psi .$ Furthermore,   $D{}
^\phi \bowtie ^\psi D$ is a double cross coproduct. We shall show
that both left coaction $\phi$ and right  coaction $\psi$ are
non-trivial.
\end {Example}

{\bf Proof.} 
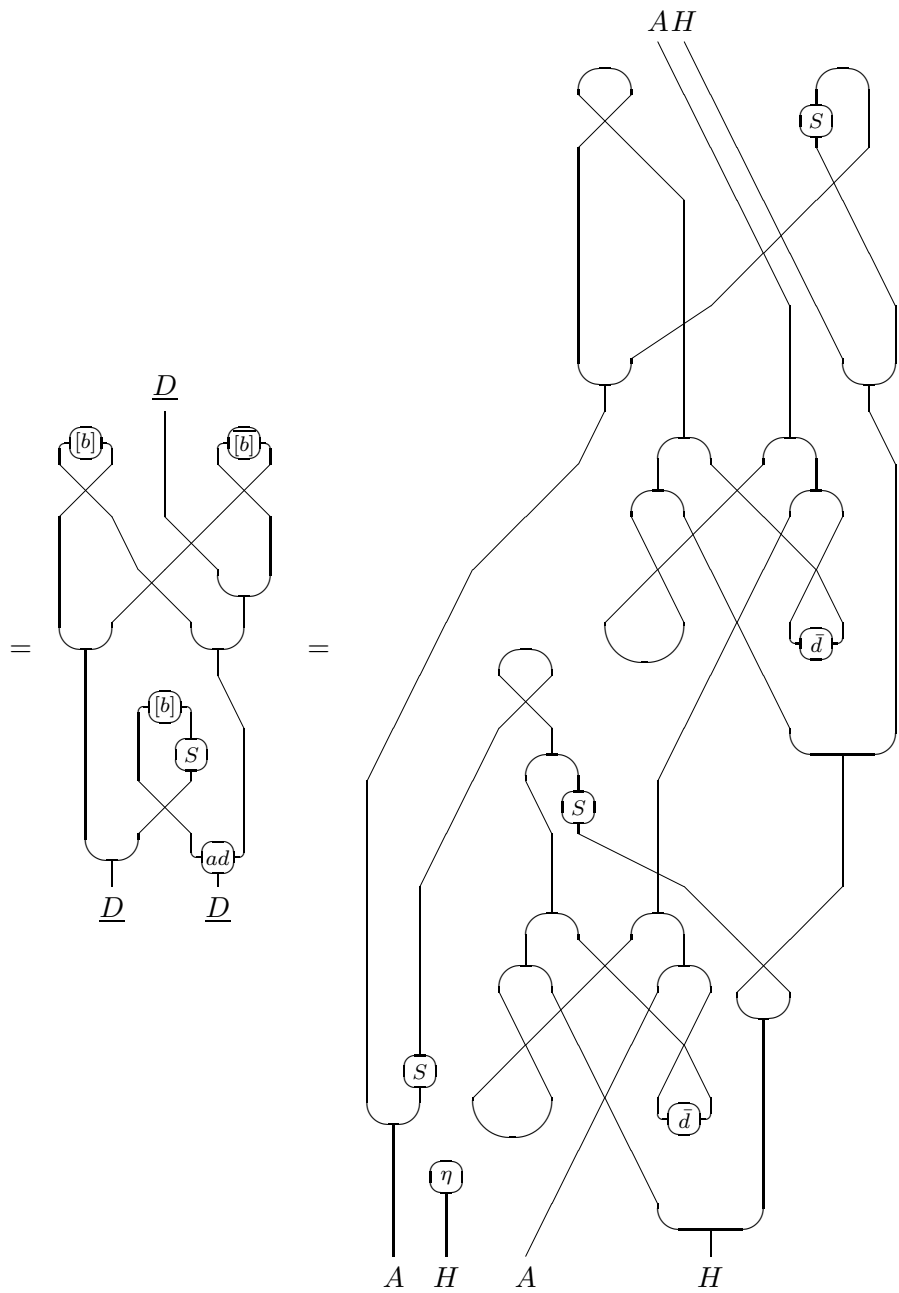
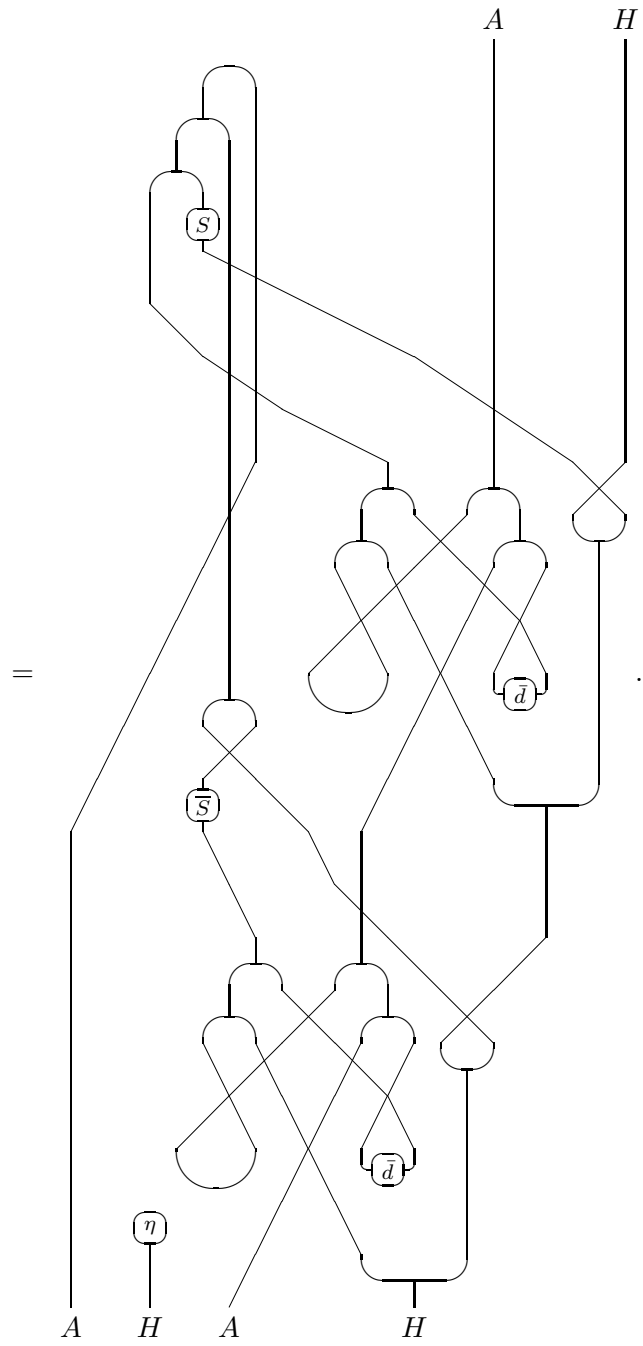
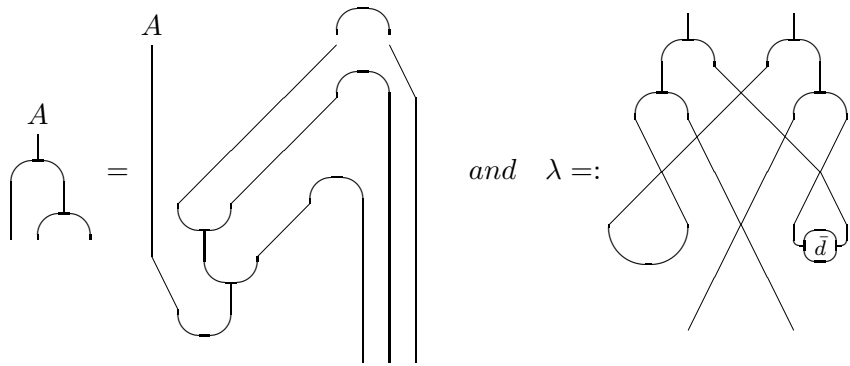
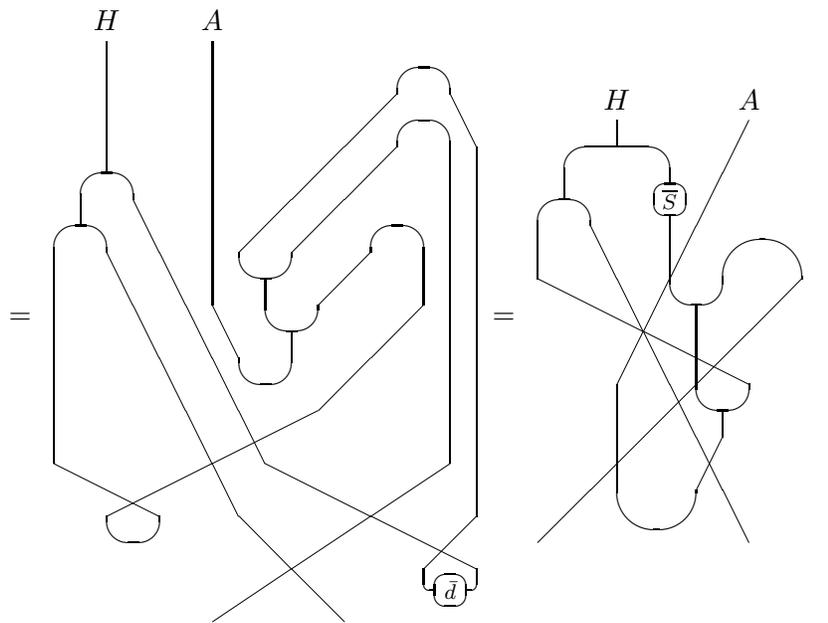
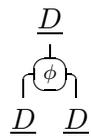
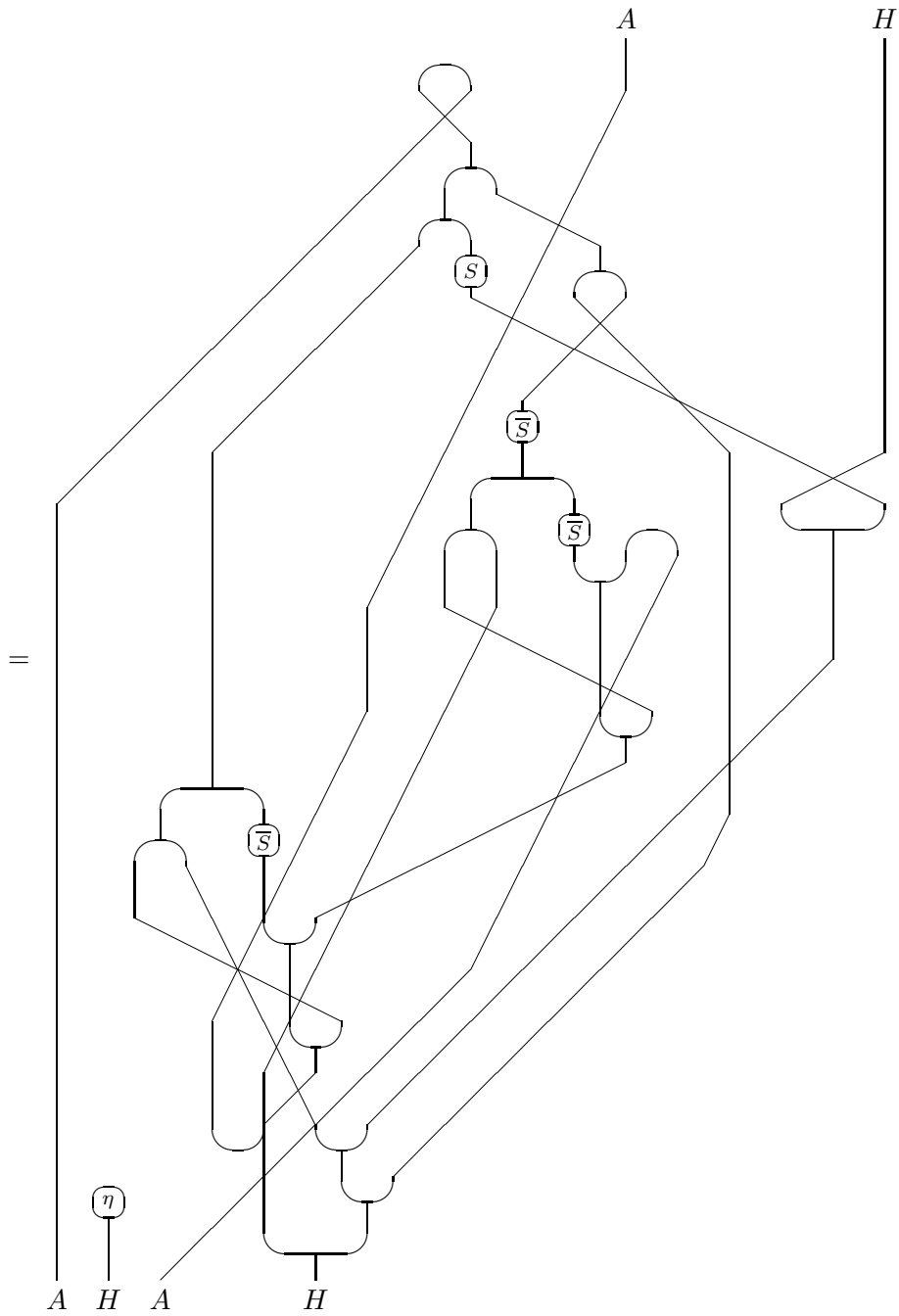
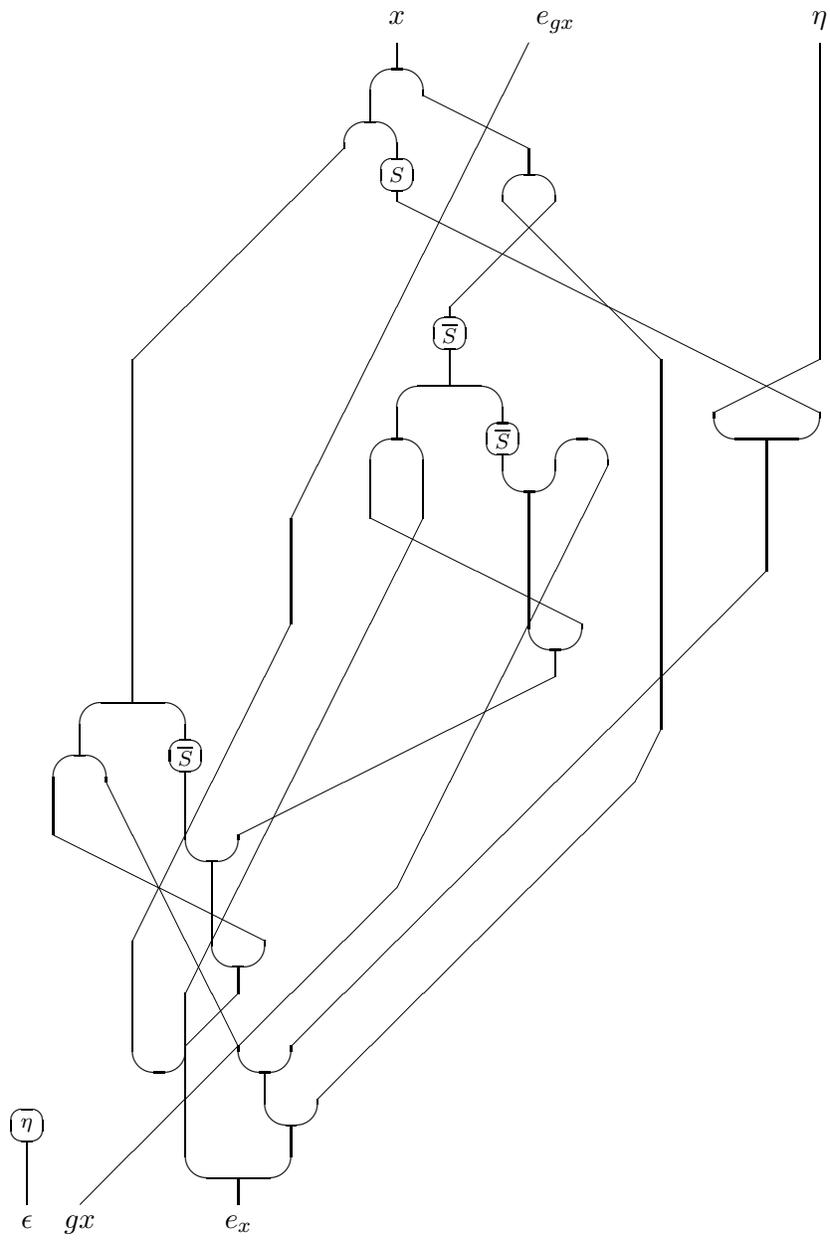
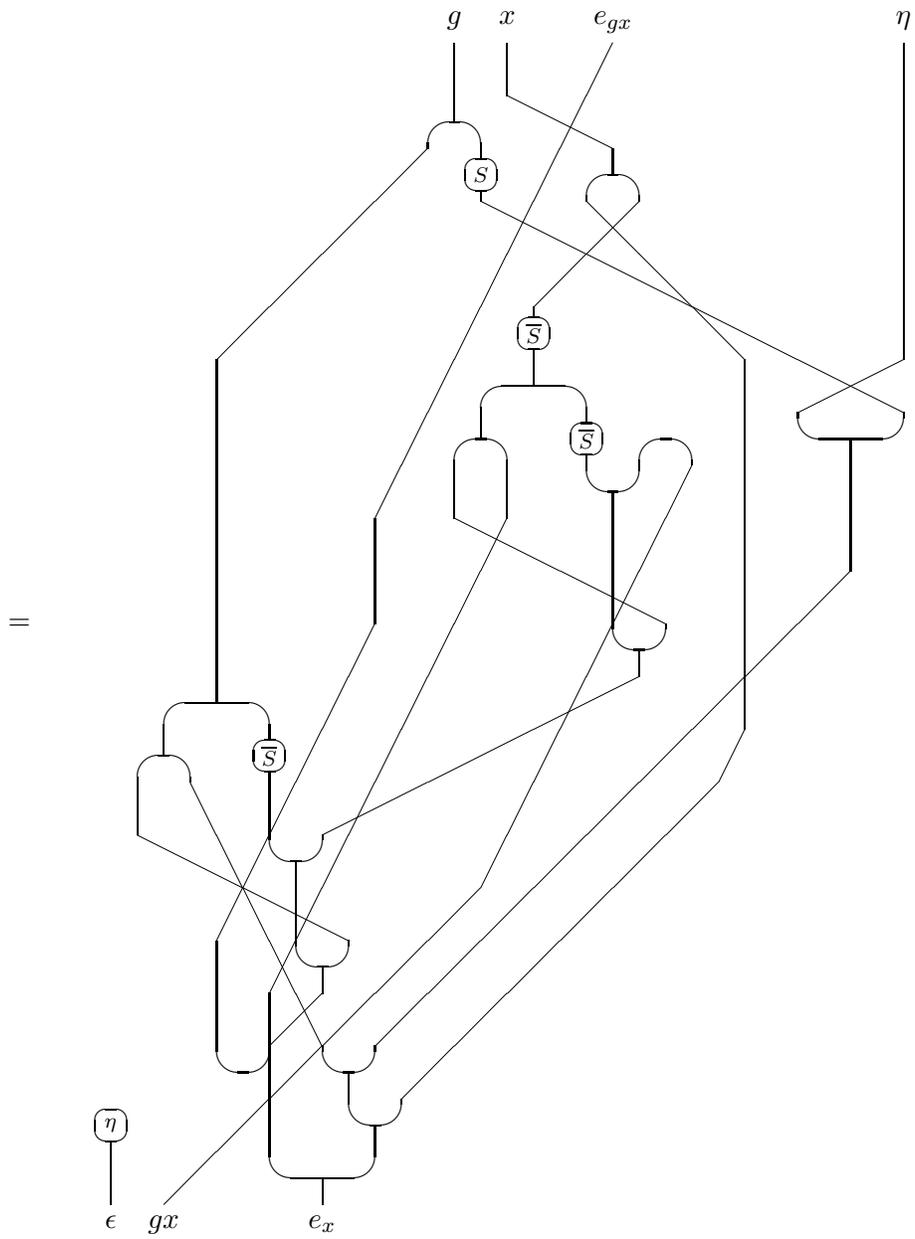
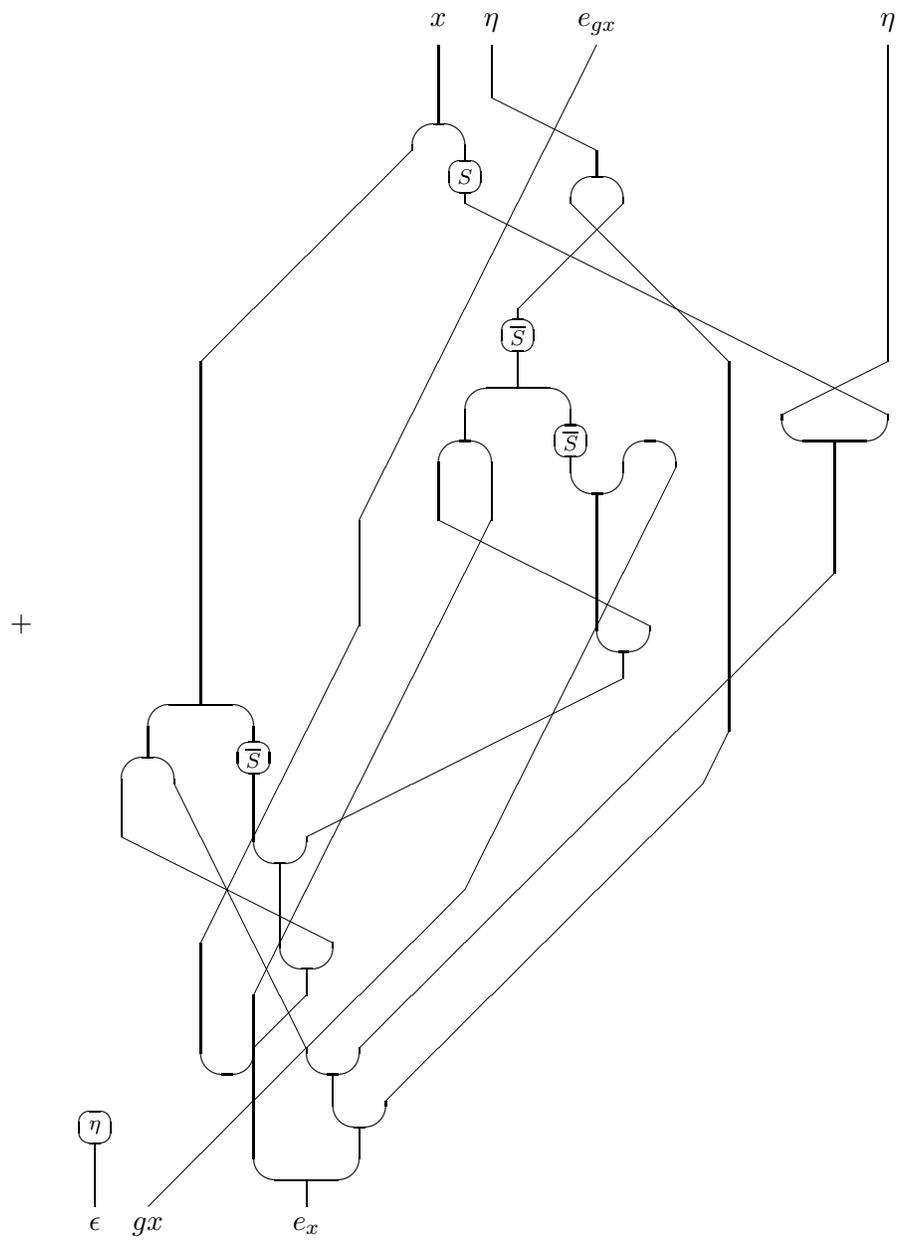
\[
\begin{tangle}
\step\object{\underline{D}}\\
\td \phi \\
\object{\underline{D}}\Step\object{\underline{D}}
 \end{tangle}\]
\[ \ \ =\ \
   \begin{tangle}
  \Step\Step\object{\underline{D}}\\
  \ro {[b]}\Step\id\Step
  \ro {\overline { [b]}}\\
  \XX\Step\id\Step\XX\\
  \id\Step\d\step\XX\Step\id\\
  \id\Step\sw2\nw2\step\step\cu\\
 \cu\Step\step\cu\\
\step\id\Step \ro {[b]}\step\d\\
\step\id\Step\id\Step\S\Step\id\\
\step\id\Step\XX\Step\id\\
\step\cu\Step\tu {ad}\\
\Step\object{\underline{D}}\Step\Step\object{\underline{D}}
   \end{tangle}
    \ \ \ = \ \ \
\begin{tangle}
\Step\Step\Step\Step\step\Step\object{A}\step\object{H}\\
\Step\Step\Step\Step\coev  \step\d\d
\Step\step\step\coev \\
\Step\Step\Step\Step\XX\Step\d\d\Step\step\S\Step\id\\
\Step\Step\Step\Step\id\Step\nw2\step\step\d\d\Step\nw1\ne2\\
\Step\Step\Step\Step\id\Step\Step\id\Step\d\nw1\ne2\step\d\\
\Step\Step\Step\Step\id\Step\Step\id\Step\ne2\d\d\Step\d\\
\Step\Step\Step\Step\id\step\step\sw3\step\id\Step\Step\id\step\d\Step\id\\
\Step\Step\Step\Step\cu\Step\id\Step\Step\id\Step\cu\\
\Step\Step\Step\Step\dd\step\step\cd\Step\cd\step\step\d\\
\Step\step\Step\Step\ne2\step\step\cd\step\XX\step\cd\Step\id\\
\Step\step\Step\ne2\Step\step\step\d\sw2\se1\step\step \sw1\se2\dd\Step\id\\
\Step\step\dd\step\Step\Step\step\ne2\nw1\step\nw1\ne1\step\ne1\nw1 \Step\id\\

\Step\dd\Step \coev\step\step\Ev\Step\step\dd\d\step\coro {\bar d} \step[2]\id\\
\step\dd\Step\step\XX\step\step\step\Step\dd\Step\d\Step\Step\id\\

\dd\step\Step\dd\step \cd\Step\step\dd\Step\Step\Cu\\

\id\Step\step\dd\Step\d\step\S\Step\step\id\Step\Step\Step\step\id\\
\id\Step\dd\Step\Step\id\step\nw4\Step\id\Step\Step\Step\step\id\\
\id\Step\id\Step\Step\cd\Step\cd\nw2\Step\Step\ne2\\
\id\Step\id\Step\step\cd\step\XX\step\cd\step\XX\\
\id\Step\id\Step\step\d\sw2\se1\step\step \sw1\se2\dd\step\cu\\
\id\Step\S\step\step\step\ne2\nw1\step\nw1\ne1\step\ne1\nw1\Step\id\\

\cu\Step\Ev\Step\step\dd\d\step \coro {\bar d}\Step\id\\
\step\id\Step\Q \eta\Step\Step\dd\Step\d\Step\Step\id\\
\step\id\Step\id\Step\step\dd\Step\Step\Cu\\
\step\object{A}\Step\object{H}\Step\step\object{A}\Step\Step\Step\step\object{H}
\end{tangle}
\]
\[
 \ \ \ = \ \ \
\begin{tangle}
\Step\step\Step\Step\Step\Step\Step\step\Step\object{A}\step\Step\Step\object{H}\\
\Step\step\Step\coev \Step\Step\step\Step\Step\id\Step\Step\step\id\\
\Step\step\step\cd \step\id\Step\Step\Step\Step\step\id\Step\Step\step\id\\
\Step\step \cd\step\id\step\id\Step\Step\Step\Step\step\id\Step\Step\step\id\\\
\Step\step\id\step\step\S\step\id\step\id\Step\Step\Step\Step\step\id\Step\Step\step\id\\\
\Step\step\id\Step\nw4\id\step\id\Step\Step\Step\Step\step\id\Step\Step\step\id\\
\Step\step\nw2\Step\id\step\id\Step\nw4\Step\Step\Step\id\Step\step\Step\id\\
\Step\step\Step\nw3\id\step\id\Step\Step\Step\nw3\Step\id\Step\Step\step\id\\
\Step\step\Step\step\id\step\id\step\nw4\Step\Step\Step\step\id\nw3\Step\Step\id\\

\Step\Step\Step\id\dd\Step\step\step\cd\Step\cd\step\step\XX\\
\Step\step\Step\dd\id\Step\step\step\cd\step\XX\step\cd\step\cu\\
\Step\Step\dd\step\id\Step\step\step\d\sw2\se1\step\step \sw1\se2\dd\Step\id\\
\Step\step\dd\step\step\id\step\Step\step\ne2\nw1\step\nw1\ne1\step\ne1\nw1\Step\id\\
\Step\dd\Step  \cd\step\step\Ev\Step\step\dd\d\step\coro {\bar d}\Step\id\\
\step\dd\Step\step\XX\step\step\step\Step\dd\Step\d\Step\Step\id\\
\dd\Step\Step\O {\overline{S}}\Step\nw2\Step\step\dd\Step\Step \Cu\\

\id\Step\step\Step\d\Step\step\d\step\id\Step\Step\Step\step\id\\
\id\Step\Step\Step\d\step\Step\nw2\id\Step\Step\Step\step\id\\
\id\Step\Step\Step\cd\Step\cd\nw2\Step\Step\ne2\\
\id\Step\Step\step\cd\step\XX\step\cd\step\XX\\
\id\Step\Step\step\d\sw2\se1\step\step \sw1\se2\dd\step\cu\\
\id\Step\step\step\step\ne2\nw1\step\nw1\ne1\step\ne1\nw1\Step\id\\
\id\Step\Step\Ev\Step\step\dd\d\step\coro {\bar d }\Step\id\\
\id\step\Step\Q \eta\Step\Step\dd\Step\d\Step\Step\id\\
\id\step\Step\id\Step\step\dd\Step\Step\Cu\\
\object{A}\step\Step\object{H}\Step\step\object{A}\Step\Step\Step\step\object{H}
\end{tangle} \ .
\] Compute:  \[
\begin{tangle}
\step\object{A}\\
\cd\\
\id \step\cd\\
\end{tangle}
\;=\enspace
\begin{tangle}
\object{A}\Step\Step\Step\step\coev\\
\id\Step\Step\Step\ne2\coev \d\\
\id\Step\Step\ne2\step\ne2\Step\id\step\id\\
\id\Step\ne2\step\ne2\step\coev\step\id\step\id\\
\id\step\cu\Step\ne2\Step\id\step\id\step\id\\
\d\step\cu\Step\Step\id\step\id\step\id\\
\step\ev\Step\Step\step\id\step\id\step\id\\
\end{tangle}
\Step and \ \ \  \lambda =:
\begin{tangle}
\step\step \cd\Step\cd\\
\step\cd\step
\XX\step\cd\\
\step\d\sw2\se1\step\step \sw1\se2\dd\\
\step\ne2\nw1\step\nw1\ne1\step\ne1\nw1\\
\Ev\Step\step\dd\d\step\coro {\bar d }\Step\\
\step\Step\dd\Step\d\Step\\
\end{tangle}
\Step
\]\[
\;=\enspace
\begin{tangle}
\Step\object{H}\Step\Step\object{A}\\
\Step\id\Step\Step\id\Step\Step\Step\step\coev\\
\Step\id\Step\Step\id\Step\Step\Step\ne2\coev\d\\
\step\cd\step\Step\id\Step\Step\ne2\step\ne2\Step\id\step\id\\
\cd\step\d\Step\id\Step\ne2\step\ne2\step\coev\step\id\step\id\\
\id\Step\d\step\d\step\id\step\cu\Step\ne2\Step\id\step\id\step\id\\
\id\Step\step\d\step\d\d\step\cu\Step\step\ne2\step\id\step\id\\
\id\Step\Step\d\step\d\ev\Step\ne2\step\Step\id\step\id\\
\id\Step\Step\step\d\sw4\se1\Step\Step\Step\step\id\step\id\\
\nw4\step\sw4\Step\step\d\step\nw4\Step\step\sw3\Step\step\id\\
\Step\ev\Step\Step\se2\sw3\Step\nw4\step\sw2\\
\Step\Step\Step\Step\ne3\nw2\Step\Step\coro {\bar d}\\
\end{tangle}
\;=\enspace
\begin{tangle}
\Step\step\object{H}\Step\Step\step\object{A}\\
\step\Cd\Step\dd\\
\cd\Step\step\O {\overline{S}}\step\dd\\
\id\Step\d\Step\id\dd\step\Coev\\
\nw4\Step\nw1\ne1\cu\Step\ne2\\
\Step\step\dd\se1\step\id\sw2\se4\\
\Step\step\id\step\sw2\se1\cu\\
\Step\sw2\id\Step\step\X\\
\sw2\Step\Ev \Step\step\d
\end{tangle} \ \ \ .
\] Also,
\[
\begin{tangle}
\step\object{\underline{D}}\\
\td \phi\\
\object{\underline{D}}\Step\object{\underline{D}}\\
\end{tangle} \]\[
= \ \
\begin{tangle}
\Step\Step\Step\Step\Step\Step\Step\Step\Step\Step\Step
\object{A}\Step\Step\Step\Step\Step\object{H}\\
\Step\Step\Step\Step\Step\Step\Step\coev\Step\Step\Step
\id\Step\Step\Step\Step\Step\id\\  
\Step\Step\Step\Step\Step\Step\Step\XX
\Step\Step\step\dd\Step\Step\Step\Step\Step\id\\   
\Step\Step\Step\Step\Step\Step\step\ne2\step\cd\Step\step\dd
\Step\Step\Step\Step\Step\step\id\\   
\Step\Step\Step\Step\Step\step\ne2\Step\cd
\step\nw4\step\sw1\Step\Step\Step\Step\Step\Step\id\\ 
\Step\Step\Step\Step\step\ne2\Step\step\ne2\Step\S\Step\dd\step\cd
\Step\Step\Step\Step\Step\id\\   
\Step\Step\Step\step\ne2\step\Step\ne2\Step\Step\nw4\ne1\Step\XX
\Step\Step\Step\Step\Step\id\\  
\Step\Step\step\ne2\step\Step\ne2\Step\Step\Step\dd\Step\ne2\nw4\step\nw2
\Step\Step\Step\Step\step\id\\    
\Step\step\ne2\step\Step\ne2\Step\Step\Step\step\dd\Step\O
{\overline{S}}\Step\Step\step\Step\se2\step\se4\Step\Step\id\\   
\step\ne2\Step\Step\id\Step\Step\Step\Step\dd\step\Cd\Step\Step\Step
\id\Step\nw4\Step\ne4\\    
\id\Step\Step\Step\id\Step\Step\Step\step\dd\step\cd \Step\step\O
{\overline{S}}\Step\coev\Step\id\Step\Cu\\   
\id\Step\Step\Step\id\Step\Step\Step\dd
\Step\id\Step\id\Step\step\cu\step\dd\Step\id\Step\Step\id\\   
\id\Step\Step\Step\id\Step\Step\Step\id\Step\step\nw4\sw1\Step\Step
\id\step\dd\Step\step\id\Step\Step\id\\      
\id\Step\Step\Step\id\Step\Step\Step\id\Step\step\dd
\Step\Step\step\id\sw1\se4\Step\step\id\Step\step\ne2\\  
\id\Step\Step\Step\id\Step\Step\step\dd
\Step\dd\Step\step\Step\dd\cu\Step\step\id\step\ne2\\  
\id\Step\Step\Cd\Step\dd\Step\dd\Step\Step
\sw4\sw1\Step\Step\step\ne2\id\\           
\id\step\Step\cd\step\Step\O
{\overline{S}}\step\dd\Step\dd\Step\Step\ne4\dd
\Step\Step\ne2\step\dd\\  
\id\Step\step\id\Step\d\Step\id\dd\step\sw4\sw1\Step\step\Step
\dd\Step\step\ne2\Step\ne2\\   
\id\Step\step\nw4\Step\nw1\ne1\cu\dd
\Step\Step\step\dd\Step\ne2\Step\ne2\\   
\id\Step\Step\Step\dd\d\step\id\sw1\se4
\Step\Step\ne2\step\ne2\step\step\ne2\\     
\id\Step\Step\Step\id\Step\X\cu\Step\ne2\step\ne2\step\step\ne2\\ 
\id\Step\Step\Step\id\Step\id\sw2\nw1
\step\ne2\step\ne2\step\step\ne2\\  
\id\Step\Step\Step\ev\id\step\ne2\cu\Step\ne2\\  
\id\Step \Q {\eta} \step\Step\Step\ne2\id\Step\step\cu\\  
\id\Step\id\Step\step\ne2\Step\Cu\\  
\object{A}\Step\object{H}\Step\object{A}\Step\Step\Step\object{H}\\
\end{tangle}
\] and \[
\begin{tangle}
\Step\object{e_{gx}\otimes \eta}\\
\td \phi \\
\id\Step\id\\
\object{x\otimes \epsilon}\Step\step\object{gx\otimes e_{x}}\\
\end{tangle}
\Step\;=\enspace
\]\[
\begin{tangle}
\Step\Step\Step\Step\Step\Step\Step\Step\object{x}\Step\Step\Step
\object{e_{gx}}\Step\Step\Step\Step\Step\object{\eta}\\
\Step\Step\Step\Step\Step\Step\step\step\step\cd\Step\step\dd
\Step\Step\Step\Step\Step\step\id\\   
\Step\Step\Step\Step\Step\step\step\Step\cd
\step\nw4\step\sw1\Step\Step\Step\Step\Step\Step\id\\ 
\Step\Step\Step\Step\step\step\Step\step\ne2\Step\S\Step\dd\step\cd
\Step\Step\Step\Step\Step\id\\   
\Step\Step\Step\step\step\step\Step\ne2\Step\Step\nw4\ne1\Step\XX
\Step\Step\Step\Step\Step\id\\  
\Step\Step\step\step\step\Step\ne2\Step\Step\Step\dd\Step\ne2\nw4\step\nw2
\Step\Step\Step\Step\step\id\\    
\Step\step\step\step\Step\ne2\Step\Step\Step\step\dd\Step\O
{\overline{S}}\Step\Step\step\Step\se2\step\se4\Step\Step\id\\   
\step\step\Step\Step\id\Step\Step\Step\Step\dd\step\Cd
\Step\Step\Step\id\Step\nw4\Step\ne4\\    
\Step\Step\Step\id\Step\Step\Step\step\dd\step\cd \Step\step\O
{\overline{S}}\Step\coev\Step\id\Step\Cu\\   
\Step\Step\Step\id\Step\Step\Step\dd
\Step\id\Step\id\Step\step\cu\step\dd\Step\id\Step\Step\id\\   
\Step\Step\Step\id\Step\Step\Step\id\Step\step\nw4\sw1
\Step\Step\id\step\dd\Step\step\id\Step\Step\id\\      
\Step\Step\Step\id\Step\Step\Step\id\Step\step\dd
\Step\Step\step\id\sw1\se4\Step\step\id\Step\step\ne2\\  
\Step\Step\Step\id\Step\Step\step\dd
\Step\dd\Step\step\Step\dd\cu\Step\step\id\step\ne2\\  
\Step\Step\Cd\Step\dd\Step\dd\Step\Step
\sw4\sw1\Step\Step\step\ne2\id\\           
\step\Step\cd\step\Step\O
{\overline{S}}\step\dd\Step\dd\Step\Step\ne4\dd
\Step\Step\ne2\step\dd\\  
\Step\step\id\Step\d\Step\id\dd\step\sw4\sw1\Step\step\Step
\dd\Step\step\ne2\Step\ne2\\   
\Step\step\nw4\Step\nw1\ne1\cu\dd
\Step\Step\step\dd\Step\ne2\Step\ne2\\   
\Step\Step\Step\dd\d\step\id\sw1\se4
\Step\Step\ne2\step\ne2\step\step\ne2\\     
\Step\Step\Step\id\Step\X\cu\Step\ne2\step\ne2\step\step\ne2\\ 
\Step\Step\Step\id\Step\id\sw2\nw1
\step\ne2\step\ne2\step\step\ne2\\  
\Step\Step\Step\ev\id\step\ne2\cu\Step\ne2\\  
\Step \Q {\eta} \step\Step\Step\ne2\id\Step\step\cu\\  
\Step\id\Step\step\ne2\Step\Cu\\  
\Step\object{\epsilon}\Step\object{gx}\Step\Step\Step\object{e_{x}}\\
\end{tangle}
\]\[
= \ \
\begin{tangle}
\Step\Step\Step\Step\Step\Step\Step\step\object{g}\Step\object{x}\Step\Step
\object{e_{gx}}\Step\Step\Step\Step\Step\step\object{\eta}\\
\Step\Step\Step\Step\Step\Step\step\step\step\id\step\step\id\Step\step\dd
\Step\Step\Step\Step\Step\step\id\\   
\Step\Step\Step\Step\Step\step\step\Step\cd
\step\nw4\step\sw1\Step\Step\Step\Step\Step\Step\id\\ 
\Step\Step\Step\Step\step\step\Step\step\ne2\Step\S\Step\dd\step\cd
\Step\Step\Step\Step\Step\id\\   
\Step\Step\Step\step\step\step\Step\ne2\Step\Step\nw4\ne1\Step\XX
\Step\Step\Step\Step\Step\id\\  
\Step\Step\step\step\step\Step\ne2\Step\Step\Step\dd\Step\ne2\nw4\step\nw2
\Step\Step\Step\Step\step\id\\    
\Step\step\step\step\Step\ne2\Step\Step\Step\step\dd\Step\O
{\overline{S}}\Step\Step\step\Step\se2\step\se4\Step\Step\id\\   
\step\step\Step\Step\id\Step\Step\Step\Step\dd\step\Cd
\Step\Step\Step\id\Step\nw4\Step\ne4\\    
\Step\Step\Step\id\Step\Step\Step\step\dd\step\cd \Step\step\O
{\overline{S}}\Step\coev\Step\id\Step\Cu\\   
\Step\Step\Step\id\Step\Step\Step\dd
\Step\id\Step\id\Step\step\cu\step\dd\Step\id\Step\Step\id\\   
\Step\Step\Step\id\Step\Step\Step\id\Step\step\nw4\sw1
\Step\Step\id\step\dd\Step\step\id\Step\Step\id\\      
\Step\Step\Step\id\Step\Step\Step\id\Step\step\dd
\Step\Step\step\id\sw1\se4\Step\step\id\Step\step\ne2\\  
\Step\Step\Step\id\Step\Step\step\dd
\Step\dd\Step\step\Step\dd\cu\Step\step\id\step\ne2\\  
\Step\Step\Cd\Step\dd\Step\dd\Step\Step
\sw4\sw1\Step\Step\step\ne2\id\\           
\step\Step\cd\step\Step\O
{\overline{S}}\step\dd\Step\dd\Step\Step\ne4\dd
\Step\Step\ne2\step\dd\\  
\Step\step\id\Step\d\Step\id\dd\step\sw4\sw1\Step\step\Step
\dd\Step\step\ne2\Step\ne2\\   
\Step\step\nw4\Step\nw1\ne1\cu\dd
\Step\Step\step\dd\Step\ne2\Step\ne2\\   
\Step\Step\Step\dd\d\step\id\sw1\se4
\Step\Step\ne2\step\ne2\step\step\ne2\\     
\Step\Step\Step\id\Step\X\cu\Step\ne2\step\ne2\step\step\ne2\\ 
\Step\Step\Step\id\Step\id\sw2\nw1
\step\ne2\step\ne2\step\step\ne2\\  
\Step\Step\Step\ev\id\step\ne2\cu\Step\ne2\\  
\Step \Q {\eta} \step\Step\Step\ne2\id\Step\step\cu\\  
\Step\id\Step\step\ne2\Step\Cu\\  
\Step\object{\epsilon}\Step\object{gx}\Step\Step\Step\object{e_{x}}\\
\end{tangle}\]\[
+ \
\begin{tangle}
\Step\Step\Step\Step\Step\Step\Step\step\object{x}\Step\object{\eta}\Step\Step
\object{e_{gx}}\Step\Step\Step\Step\Step\step\object{\eta}\\
\Step\Step\Step\Step\Step\Step\step\step\step\id\step\step\id\Step\step\dd
\Step\Step\Step\Step\Step\step\id\\   
\Step\Step\Step\Step\Step\step\step\Step\cd
\step\nw4\step\sw1\Step\Step\Step\Step\Step\Step\id\\ 
\Step\Step\Step\Step\step\step\Step\step\ne2\Step\S\Step\dd\step\cd
\Step\Step\Step\Step\Step\id\\   
\Step\Step\Step\step\step\step\Step\ne2\Step\Step\nw4\ne1\Step\XX
\Step\Step\Step\Step\Step\id\\  
\Step\Step\step\step\step\Step\ne2\Step\Step\Step\dd\Step\ne2\nw4\step\nw2
\Step\Step\Step\Step\step\id\\    
\Step\step\step\step\Step\ne2\Step\Step\Step\step\dd\Step\O
{\overline{S}}\Step\Step\step\Step\se2\step\se4\Step\Step\id\\   
\step\step\Step\Step\id\Step\Step\Step\Step\dd\step\Cd
\Step\Step\Step\id\Step\nw4\Step\ne4\\    
\Step\Step\Step\id\Step\Step\Step\step\dd\step\cd \Step\step\O
{\overline{S}}\Step\coev\Step\id\Step\Cu\\   
\Step\Step\Step\id\Step\Step\Step\dd
\Step\id\Step\id\Step\step\cu\step\dd\Step\id\Step\Step\id\\   
\Step\Step\Step\id\Step\Step\Step\id\Step\step\nw4\sw1
\Step\Step\id\step\dd\Step\step\id\Step\Step\id\\      
\Step\Step\Step\id\Step\Step\Step\id\Step\step\dd
\Step\Step\step\id\sw1\se4\Step\step\id\Step\step\ne2\\  
\Step\Step\Step\id\Step\Step\step\dd
\Step\dd\Step\step\Step\dd\cu\Step\step\id\step\ne2\\  
\Step\Step\Cd\Step\dd\Step\dd\Step\Step
\sw4\sw1\Step\Step\step\ne2\id\\           
\step\Step\cd\step\Step\O
{\overline{S}}\step\dd\Step\dd\Step\Step\ne4\dd
\Step\Step\ne2\step\dd\\  
\Step\step\id\Step\d\Step\id\dd\step\sw4\sw1\Step\step\Step
\dd\Step\step\ne2\Step\ne2\\   
\Step\step\nw4\Step\nw1\ne1\cu\dd
\Step\Step\step\dd\Step\ne2\Step\ne2\\   
\Step\Step\Step\dd\d\step\id\sw1\se4
\Step\Step\ne2\step\ne2\step\step\ne2\\     
\Step\Step\Step\id\Step\X\cu\Step\ne2\step\ne2\step\step\ne2\\ 
\Step\Step\Step\id\Step\id\sw2\nw1
\step\ne2\step\ne2\step\step\ne2\\  
\Step\Step\Step\ev\id\step\ne2\cu\Step\ne2\\  
\Step \Q {\eta} \step\Step\Step\ne2\id\Step\step\cu\\  
\Step\id\Step\step\ne2\Step\Cu\\  
\Step\object{\epsilon}\Step\object{gx}\Step\Step\Step\object{e_{x}}\\
\end{tangle}
\]
Let $u$ and $v$ denote the first term and the second term,
respectively.
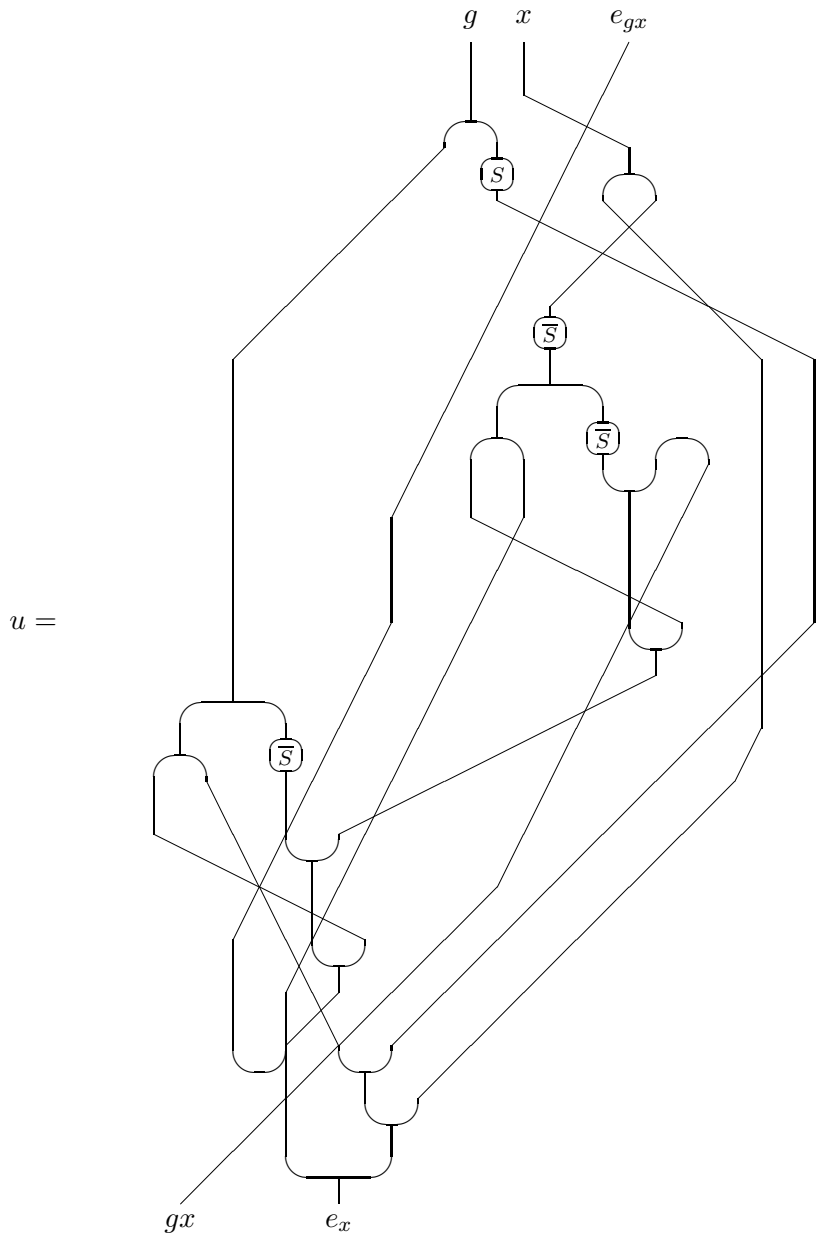
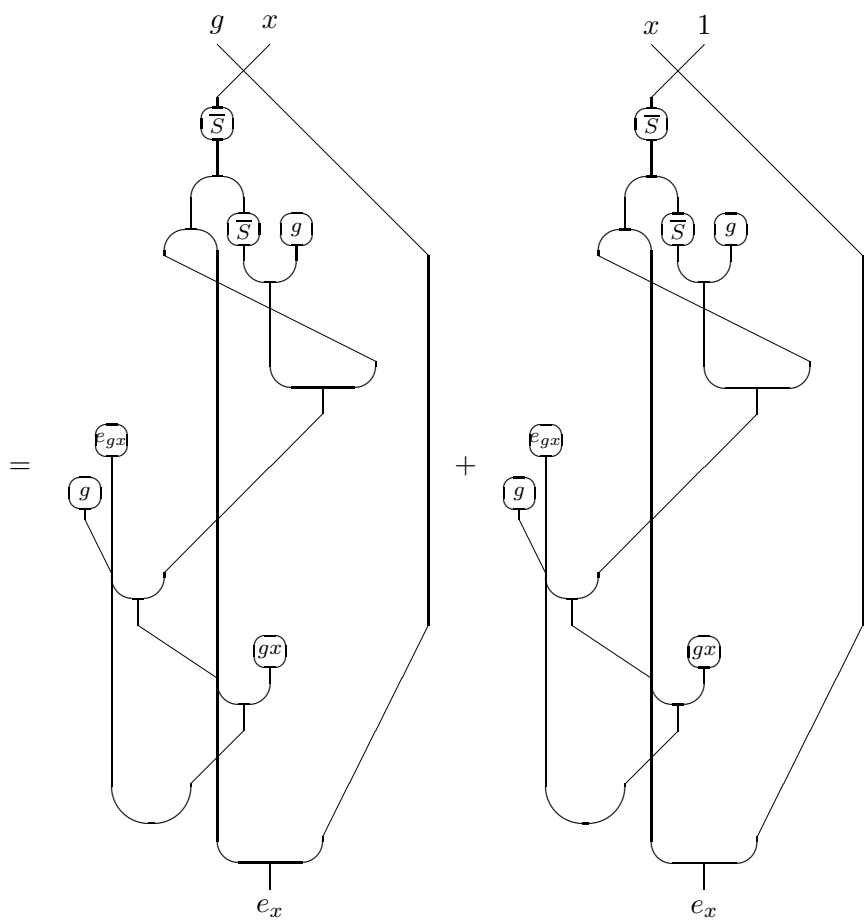
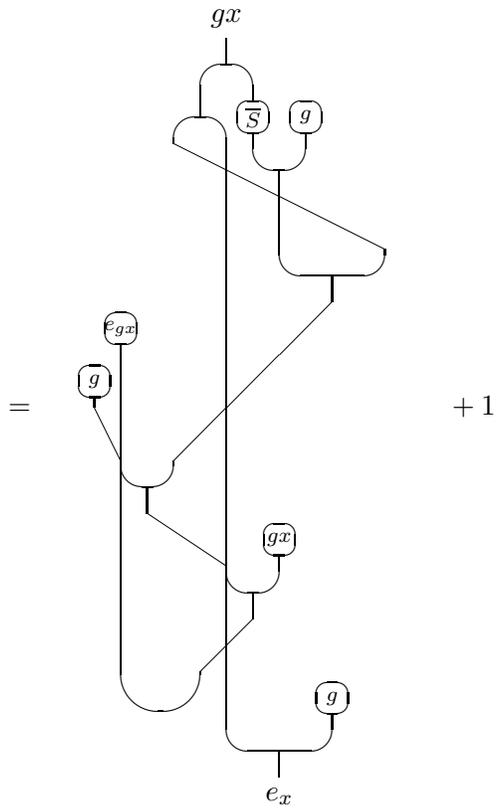
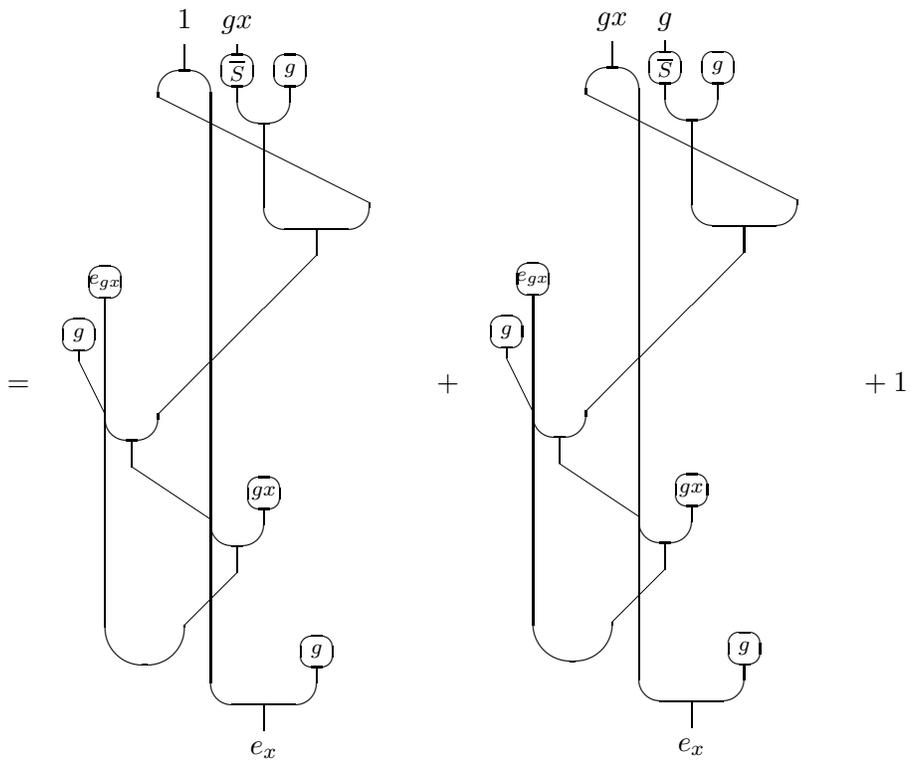
\[
u = \
\begin{tangle}
\Step\Step\Step\Step\Step\Step\Step\step\object{g}\Step\object{x}\Step\Step
\object{e_{gx}}\\
\Step\Step\Step\Step\Step\Step\step\step\step\id\step\step\id\Step\step\dd
\Step\Step\Step\Step\Step\step\\   
\Step\Step\Step\Step\Step\step\step\Step\cd
\step\nw4\step\sw1\Step\Step\Step\Step\Step\Step\\ 
\Step\Step\Step\Step\step\step\Step\step\ne2\Step\S\Step\dd\step\cd
\Step\Step\Step\Step\Step\\   
\Step\Step\Step\step\step\step\Step\ne2\Step\Step\nw4\ne1\Step\XX
\Step\Step\Step\Step\Step\\  
\Step\Step\step\step\step\Step\ne2\Step\Step\Step\dd\Step\ne2\nw4\step\nw2
\Step\Step\Step\Step\step\\    
\Step\step\step\step\Step\ne2\Step\Step\Step\step\dd\Step\O
{\overline{S}}\Step\Step\step\Step\se2\step\se4\Step\Step\\   
\step\step\Step\Step\id\Step\Step\Step\Step\dd\step\Cd
\Step\Step\Step\id\Step\id\\    
\Step\Step\Step\id\Step\Step\Step\step\dd\step\cd \Step\step\O
{\overline{S}}\Step\coev\Step\id\Step\id\\   
\Step\Step\Step\id\Step\Step\Step\dd
\Step\id\Step\id\Step\step\cu\step\dd\Step\id\Step\id\\   
\Step\Step\Step\id\Step\Step\Step\id\Step\step\nw4\sw1
\Step\Step\id\step\dd\Step\step\id\Step\id\\      
\Step\Step\Step\id\Step\Step\Step\id\Step\step\dd
\Step\Step\step\id\sw1\se4\Step\step\id\Step\id\\  
\Step\Step\Step\id\Step\Step\step\dd
\Step\dd\Step\step\Step\dd\cu\Step\step\id\step\ne2\\  
\Step\Step\Cd\Step\dd\Step\dd\Step\Step
\sw4\sw1\Step\Step\step\ne2\id\\           
\step\Step\cd\step\Step\O {\overline{S}}
\step\dd\Step\dd\Step\Step\ne4\dd
\Step\Step\ne2\step\dd\\  
\Step\step\id\Step\d\Step\id\dd\step\sw4\sw1\Step\step\Step
\dd\Step\step\ne2\Step\ne2\\   
\Step\step\nw4\Step\nw1\ne1\cu\dd
\Step\Step\step\dd\Step\ne2\Step\ne2\\   
\Step\Step\Step\dd\d\step\id\sw1\se4
\Step\Step\ne2\step\ne2\step\step\ne2\\     
\Step\Step\Step\id\Step\X\cu\Step\ne2\step\ne2\step\step\ne2\\ 
\Step\Step\Step\id\Step\id\sw2\nw1
\step\ne2\step\ne2\step\step\ne2\\  
\Step\Step\Step\ev\id\step\ne2\cu\Step\ne2\\  
\Step\step\Step\Step\ne2\id\Step\step\cu\\  
\Step\Step\step\ne2\Step\Cu\\  
\Step\Step\object{gx}\Step\Step\Step\object{e_{x}}\\
\end{tangle}
\]\[ = \ \
\begin{tangle}
\Step\Step\Step\object{g}\Step\object{x}\\
\Step\Step\Step\XX\\  
\Step\Step\Step\O {\overline{S}}\Step\nw2\\    
\Step\Step\step\cd\Step\step\nw2\\  
\Step\Step\cd\step\O {\overline{S}}\Step \Q g\Step\step\nw2\\  

\Step\Step\nw4\step\id\step\cu\Step\Step\step\id\\ 
\Step\Step\Step\id\Step\id\nw4\Step\Step\step\id\\   
\Step\Step\Step\id\Step\Cu\Step\id\\   
\Step \Q{e_{gx }}\Step\Step\id\Step\step\ne2\step\Step\step\id\\  

\step \Q g\step\id\Step\Step\id\step\ne2\Step\Step\Step\id\\  
\step\d\id\step\Step\ne2\id\Step\Step\Step\Step\id\\  
\Step\id\cu\Step\id\Step\Step\Step\Step\id\\  
\Step\id\step\nw3\Step\id\Step \Q {gx}\Step\Step\step\dd\\  
\Step\id\Step\Step\id\cu\Step\Step\dd\\  
\Step\id\Step\step\sw2\id\Step\Step\step\dd\\  
\Step\Ev\Step\step\id\Step\Step\dd\\   
\Step\Step\Step\Cu\\   
\Step\Step\Step\Step\object{e_{x}}
\end{tangle} \ \
+ \
\begin{tangle}
\Step\Step\Step\object{x}\Step\object{1}\\
\Step\Step\Step\XX\\  
\Step\Step\Step\O {\overline{S}}\Step\nw2\\    
\Step\Step\step\cd\Step\step\nw2\\  
\Step\Step\cd\step\O {\overline{S}}\Step \Q g\Step\step\nw2\\  

\Step\Step\nw4\step\id\step\cu\Step\Step\step\id\\ 
\Step\Step\Step\id\Step\id\nw4\Step\Step\step\id\\   
\Step\Step\Step\id\Step\Cu\Step\id\\   
\Step \Q{e_{gx }}\Step\Step\id\Step\step\ne2\step\Step\step\id\\  

\step \Q g\step\id\Step\Step\id\step\ne2\Step\Step\Step\id\\  
\step\d\id\step\Step\ne2\id\Step\Step\Step\Step\id\\  
\Step\id\cu\Step\id\Step\Step\Step\Step\id\\  
\Step\id\step\nw3\Step\id\Step \Q {gx}\Step\Step\step\dd\\  
\Step\id\Step\Step\id\cu\Step\Step\dd\\  
\Step\id\Step\step\sw2\id\Step\Step\step\dd\\  
\Step\Ev\Step\step\id\Step\Step\dd\\   
\Step\Step\Step\Cu\\   
\Step\Step\Step\Step\object{e_{x}}
\end{tangle} \]\[
\ \ \ = \ \ \
\begin{tangle}
\Step\Step\Step\object{gx}\\
\Step\Step\step\cd\Step\step\\  
\Step\Step\cd\step\O {\overline{S}}\Step \Q g\Step\step\\  
\Step\Step\nw4\step\id\step\cu\Step\Step\step\\ 
\Step\Step\Step\id\Step\id\nw4\Step\Step\step\\   
\Step\Step\Step\id\Step\Cu\Step\\   

\Step \Q {e_{gx}}\Step\Step\id\Step\step\ne2\step\Step\step\\  
 \step \Q g\step\id\Step\Step\id\step\ne2\Step\Step\Step\\  
\step\d\id\step\Step\ne2\id\Step\Step\Step\Step\\  
\Step\id\cu\Step\id\Step\Step\Step\Step\\  
\Step\id\step\nw3\Step\id\Step \Q {gx}\Step\Step\step\\  
\Step\id\Step\Step\id\cu\Step\Step\\  
\Step\id\Step\step\sw2\id\Step\Step\step\\  
\Step\Ev\Step\step\id\Step\Step \Q g\\   
\Step\Step\Step\Cu\\   
\Step\Step\Step\Step\object{e_{x}}
\end{tangle}
\;+1\enspace
\]\[
\;=\enspace
\begin{tangle}
\Step\Step\step\object{1}\Step\object{gx}\\
\Step\Step\cd\step\O {\overline{S}}\Step \Q g \Step\step\\  
\Step\Step\nw4\step\id\step\cu\Step\Step\step\\ 
\Step\Step\Step\id\Step\id\nw4\Step\Step\step\\   
\Step\Step\Step\id\Step\Cu\Step\\   

\Step \Q {e_{gx}}\Step\Step\id\Step\step\ne2\step\Step\step\\  

\step \Q g\step\id\Step\Step\id\step\ne2\Step\Step\Step\\  
\step\d\id\step\Step\ne2\id\Step\Step\Step\Step\\  
\Step\id\cu\Step\id\Step\Step\Step\Step\\  
\Step\id\step\nw3\Step\id\Step \Q {gx}\Step\Step\step\\  
\Step\id\Step\Step\id\cu\Step\Step\\  
\Step\id\Step\step\sw2\id\Step\Step\step\\  
\Step\Ev\Step\step\id\Step\Step \Q g\\   
\Step\Step\Step\Cu\\   
\Step\Step\Step\Step\object{e_{x}}
\end{tangle}
\;+\enspace
\begin{tangle}
\Step\Step \step\object{gx}\Step\object{g}\\

\Step\Step\cd\step\O {\overline{S}}\Step \Q g \Step\step\\  
\Step\Step\nw4\step\id\step\cu\Step\Step\step\\ 
\Step\Step\Step\id\Step\id\nw4\Step\Step\step\\   
\Step\Step\Step\id\Step\Cu\Step\\   

\Step \Q {e_{gx}}\Step\Step\id\Step\step\ne2\step\Step\step\\  

\step \Q g\step\id\Step\Step\id\step\ne2\Step\Step\Step\\  
\step\d\id\step\Step\ne2\id\Step\Step\Step\Step\\  
\Step\id\cu\Step\id\Step\Step\Step\Step\\  
\Step\id\step\nw3\Step\id\Step \Q {gx}\Step\Step\step\\  
\Step\id\Step\Step\id\cu\Step\Step\\  
\Step\id\Step\step\sw2\id\Step\Step\step\\  
\Step\Ev\Step\step\id\Step\Step \Q g\\   
\Step\Step\Step\Cu\\   
\Step\Step\Step\Step\object{e_{x}}
\end{tangle}
\;+1\enspace
\]
\[\begin{tangle}
\Step\Step\step\Step\object{g}\\
\Step\Step\Step\step\O {\overline{S}}\Step \Q g\Step\step\\  

\Step\Step \step [2]\step\cu\Step\Step\step\\ 

\Step\Step\Step\Step\id \step \Step\step \Q {gx}\\   

\Step\Step\Step\Q g\Step\Cu\Step\\   

\Step \Q {e_{gx}}\Step\Step\id\Step\step\ne2\step\Step\step\\  
\step \Q g \step\id\Step\Step\id\step\ne2\Step\Step\Step\\  
\step\d\id\step\Step\ne2\id\Step\Step\Step\Step\\  
\Step\id\cu\Step\id\Step\Step\Step\Step\\  
\Step\id\step\nw3\Step\id\Step \Q {gx}\Step\Step\step\\  

\Step\id\Step\Step\id\cu\Step\Step\\  

\Step\id\Step\step\sw2\id\Step\Step\step\\  
\Step\Ev\Step\step\id\Step\Step \Q g\\   

\Step\Step\Step\Cu\\   
\Step\Step\Step\Step\object{e_{x}}
\end{tangle}
\;+\enspace
\begin{tangle}
\Step\Step\step\Step\object{g}\\
\Step\Step\Step\step\O {\overline{S}}\Step \Q g\Step\step\\  

\Step\Step \step [2]\step\cu\Step\Step\step\\ 

\Step\Step\Step\Step\id \step \Step\step \Q {1}\\   

\Step\Step\Step\Q {gx}\Step\Cu\Step\\   

\Step \Q {e_{gx}}\Step\Step\id\Step\step\ne2\step\Step\step\\  
\step \Q g \step\id\Step\Step\id\step\ne2\Step\Step\Step\\  
\step\d\id\step\Step\ne2\id\Step\Step\Step\Step\\  
\Step\id\cu\Step\id\Step\Step\Step\Step\\  
\Step\id\step\nw3\Step\id\Step \Q {gx}\Step\Step\step\\  

\Step\id\Step\Step\id\cu\Step\Step\\  

\Step\id\Step\step\sw2\id\Step\Step\step\\  
\Step\Ev\Step\step\id\Step\Step \Q g\\   

\Step\Step\Step\Cu\\   
\Step\Step\Step\Step\object{e_{x}}
\end{tangle}
\;+1\enspace = 1
\] and \[
v \;=\enspace
\begin{tangle}
\Step\Step\object{x}\Step\Step\Step\step\object{e_{gx}}\\
\Step\step\cd\Step\Step\dd\\
\Step\step\id\Step\S\Step\step\dd\\
\step\Cd\d\step\dd\\
\cd\Step\step\O {\overline{S}}\step\X \step [1]\Q g\\
\id\Step\d\Step\id\dd\step\X\\
\nw4\Step\nw1\ne1\cu\step\d\\

\Step\step\dd\se1\step\id\step\se4\step\d\\
\Step\step\id\step\step\se1\cu\Step\d\\
\Step\step\id\Step\step\X\Step\Step\d\\
\step\Step\Ev\Step\step\d\Step\Step\id\\
\Step\Step\Step\Step\Cu\\
\Step\Step\Step\Step\Step\object{e_{x}}
\end{tangle}
\]\[
\;=\enspace
\begin{tangle}
\Step\step \object g\Step\object x \Step\Step\step\object{e_{gx}}\\

\Step\step\id\Step\S\Step\step\dd\\
\step\Cd\d\step\dd\\
\cd\Step\step\O {\overline{S}}\step\X \step \Q g\\
\id\Step\d\Step\id\dd\step\X\\
\nw4\Step\nw1\ne1\cu\step\d\\

\Step\step\dd\se1\step\id\step\se4\step\d\\
\Step\step\id\step\step\se1\cu\Step\d\\
\Step\step\id\Step\step\X\Step\Step\d\\
\step\Step\Ev\Step\step\d\Step\Step\id\\
\Step\Step\Step\Step\Cu\\
\Step\Step\Step\Step\Step\object{e_{x}}
\end{tangle}
\;+\enspace
\begin{tangle}
\Step\step \object x\Step\object 1 \Step\Step\step\object{e_{gx}}\\

\Step\step\id\Step\S\Step\step\dd\\
\step\Cd\d\step\dd\\
\cd\Step\step\O {\overline{S}}\step\X \step \Q g\\
\id\Step\d\Step\id\dd\step\X\\
\nw4\Step\nw1\ne1\cu\step\d\\

\Step\step\dd\se1\step\id\step\se4\step\d\\
\Step\step\id\step\step\se1\cu\Step\d\\
\Step\step\id\Step\step\X\Step\Step\d\\
\step\Step\Ev\Step\step\d\Step\Step\id\\
\Step\Step\Step\Step\Cu\\
\Step\Step\Step\Step\Step\object{e_{x}}
\end{tangle}
\]\[
\;=1+\enspace
\begin{tangle}
\step \object g\Step\step \object x\Step\object 1 \Step\object{e_{gx}}\\

\cd\Step\O {\overline{S}}\step\step\X \step \Q g\\
\id\Step\d\step\d\dd\step\X\\
\nw4\Step\nw1\ne1\cu\step\d\\

\Step\step\dd\se1\step\id\step\se4\step\d\\
\Step\step\id\step\step\se1\cu\Step\d\\
\Step\step\id\Step\step\X\Step\Step\d\\
\step\Step\Ev\Step\step\d\Step\Step\id\\
\Step\Step\Step\Step\Cu\\
\Step\Step\Step\Step\Step\object{e_{x}}
\end{tangle}
\;+\enspace\begin{tangle}
\step \object x\Step\step \object1\Step\object 1 \Step\object{e_{gx}}\\

\cd\Step\O {\overline{S}}\step\step\X \step \Q g\\
\id\Step\d\step\d\dd\step\X\\
\nw4\Step\nw1\ne1\cu\step\d\\

\Step\step\dd\se1\step\id\step\se4\step\d\\
\Step\step\id\step\step\se1\cu\Step\d\\
\Step\step\id\Step\step\X\Step\Step\d\\
\step\Step\Ev\Step\step\d\Step\Step\id\\
\Step\Step\Step\Step\Cu\\
\Step\Step\Step\Step\Step\object{e_{x}}
\end{tangle}
\]
\[
\;=1+\enspace
\begin{tangle}
\step \Step\step \object x\Step\object 1 \Step\object{e_{gx}}\\

\Q x \step [2]\Q 1 \Step\O {\overline{S}}\step\step\X \step \Q g\\
\id\Step\d\step\d\dd\step\X\\
\nw4\Step\nw1\ne1\cu\step\d\\

\Step\step\dd\se1\step\id\step\se4\step\d\\
\Step\step\id\step\step\se1\cu\Step\d\\
\Step\step\id\Step\step\X\Step\Step\d\\
\step\Step\Ev\Step\step\d\Step\Step\id\\
\Step\Step\Step\Step\Cu\\
\Step\Step\Step\Step\Step\object{e_{x}}
\end{tangle}
\;+\enspace
\begin{tangle}
\step \Step\step \object x\Step\object 1 \Step\object{e_{gx}}\\

\Q g \step [2]\Q x \Step\O {\overline{S}}\step\step\X \step \Q g\\
\id\Step\d\step\d\dd\step\X\\
\nw4\Step\nw1\ne1\cu\step\d\\

\Step\step\dd\se1\step\id\step\se4\step\d\\
\Step\step\id\step\step\se1\cu\Step\d\\
\Step\step\id\Step\step\X\Step\Step\d\\
\step\Step\Ev\Step\step\d\Step\Step\id\\
\Step\Step\Step\Step\Cu\\
\Step\Step\Step\Step\Step\object{e_{x}}
\end{tangle}
\ = \ 1.\] \[Thus\Step
\begin{tangle}
\Step\object{e_{gx}\otimes \eta}\\
\td \phi\\
\id\Step\id\\
\object{x\otimes \epsilon}\Step \step\object{gx\otimes e_{x}}\\
\end{tangle}
\Step\;=2 , \Step\Step but\Step\Step
\begin{tangle}
\object{\eta_{D}}\Step\step\step\object{e_{gx}\otimes \eta}\\
\id\Step\step\step\id\\
\object{x\otimes \epsilon}\Step\step\step\object{gx\otimes e_{x}}
\end{tangle}
 \Step\;=0.\enspace
\]
Consequently, $\phi $ is not  trivial.
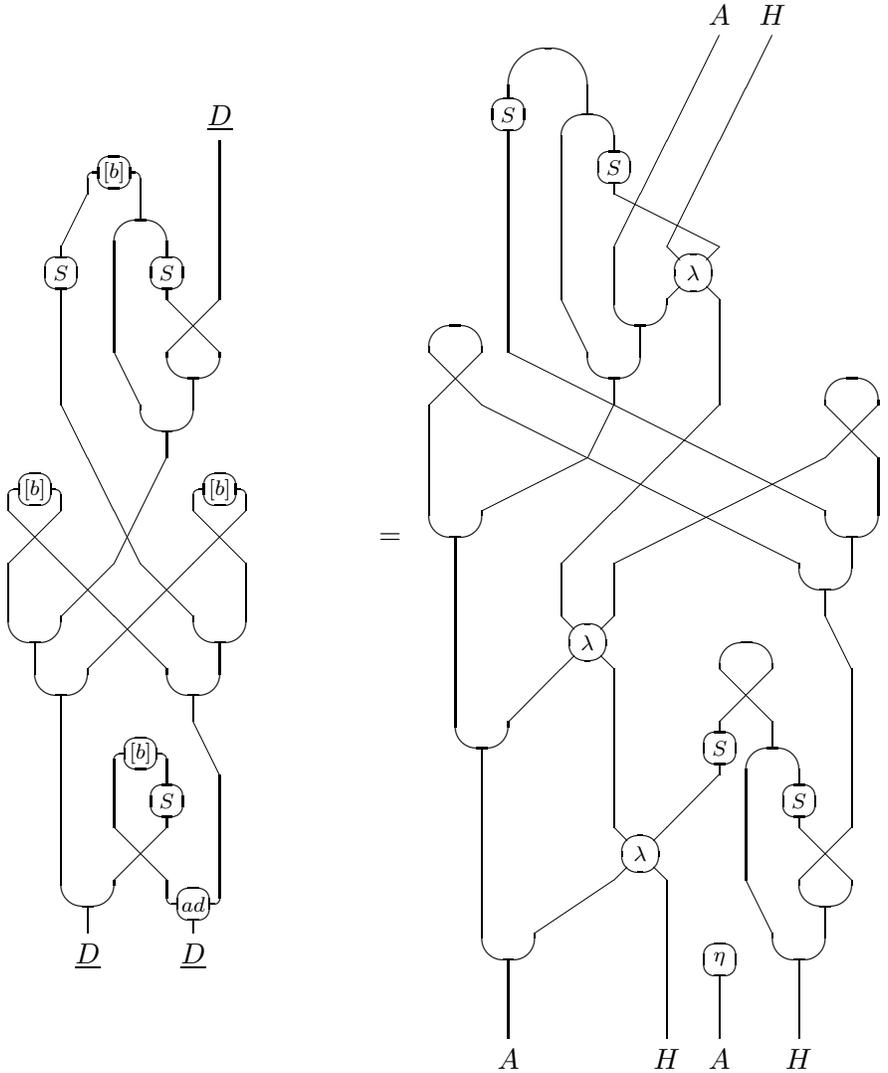
\[
\begin{tangle}
\step\object{\underline{D}}\\
\td \psi \\
\object{\underline{D}}\Step\object{\underline{D}}
 \end{tangle} \ \ = \ \ \]
\[
  \begin{tangle}
\Step\Step\Step\Step\step\object{\underline{D}}\\

\step[4] \ro {[b]}\step[3] \id\Step\Step\step\\

\Step\step\ne1 \step [1]\cd\Step\id\\
\Step\step\S \step\step\id\step\step\S\Step\id\\

\Step\step\id\step\step\id\Step\XX\\
\Step\step\id\Step\d\step\cu\\
\Step\step\d\Step\cu\\

\step\ro {[b]}\step\d\step\dd\step
 \ro {[b]}\\
 \step \XX\Step\X\Step\XX\\
 \step\id\Step\XX\step\XX\Step\id\\
 \step\cu\step\sw2\nw2\step\step\cu\\
\step\step\cu\Step\step\cu\\
\Step\step\id\Step \ro {[b]}\step\d\\
\Step\step\id\Step\id\Step\S\Step\id\\
\Step\step\id\Step\XX\Step\id\\
\Step\step\cu\Step\tu {ad}\\
\Step\Step\object{\underline{D}}\Step\Step\object{\underline{D}}
\end{tangle} \ \ = \ \
\begin{tangle}
\Step\Step\step\Step\Step\Step\object{A}\Step\object{H}\\
\Step\step\Coev\Step\Step\Step\dd\step\dd\\  
\Step\step\S\Step\cd\Step\dd\step\dd\\  
\Step\step\id\Step\id\Step\S\step\dd\step\dd\\    
\Step\step\id\Step\id\Step\ne1\step\ne1\se4\\   
\Step\step\id\Step\id\Step\id\Step\ox \lambda\\   
\coev \step\id\Step\d\step\cu\Step\id\\  
\XX\step\nw4\Step\cu\Step\step\id\Step\Step\coev\\    
\id\Step\nw4\Step\step\ne1\nw4\Step\ne2\Step\Step\XX\\   
\id\Step\Step\step\ne4\nw4\step\ne2\Step\nw4\Step\ne4\Step\id\\   
\cu\Step\Step\ne2\sw4\Step\nw4\Step\Step\cu\\    
\step\id\Step\Step\id\Step\id\Step\Step\Step\step\cu\\   
\step\id\Step\Step\ox \lambda \Step\Step\coev\Step\d\\    
\step\id\Step\step\ne2\Step\id\Step\Step\XX\Step\step\id\\    
\step\cu\Step\Step\id\Step\Step\S\step\cd\Step\id\\       
\Step\id\Step\Step\step\id\Step\step\ne2\step\id\Step\S\Step\id\\  
\Step\id\Step\Step\step\ox \lambda  \Step\step\id\Step\XX\\     
\Step\id\Step\Step\ne3\Step\id\Step\step\d\step\cu\\     
\Step\cu\Step\step\Step\id\Step\Q \eta \Step\cu\\      
\Step\step\id\Step\Step\Step\id\Step\id\Step\step\id\\   
\Step\step\object{A}\Step\Step\Step\object{H}\Step\object{A}\Step\step\object{H}
\end{tangle} \ .
\]
For convenience, we omit $g$ in the following diagrams:
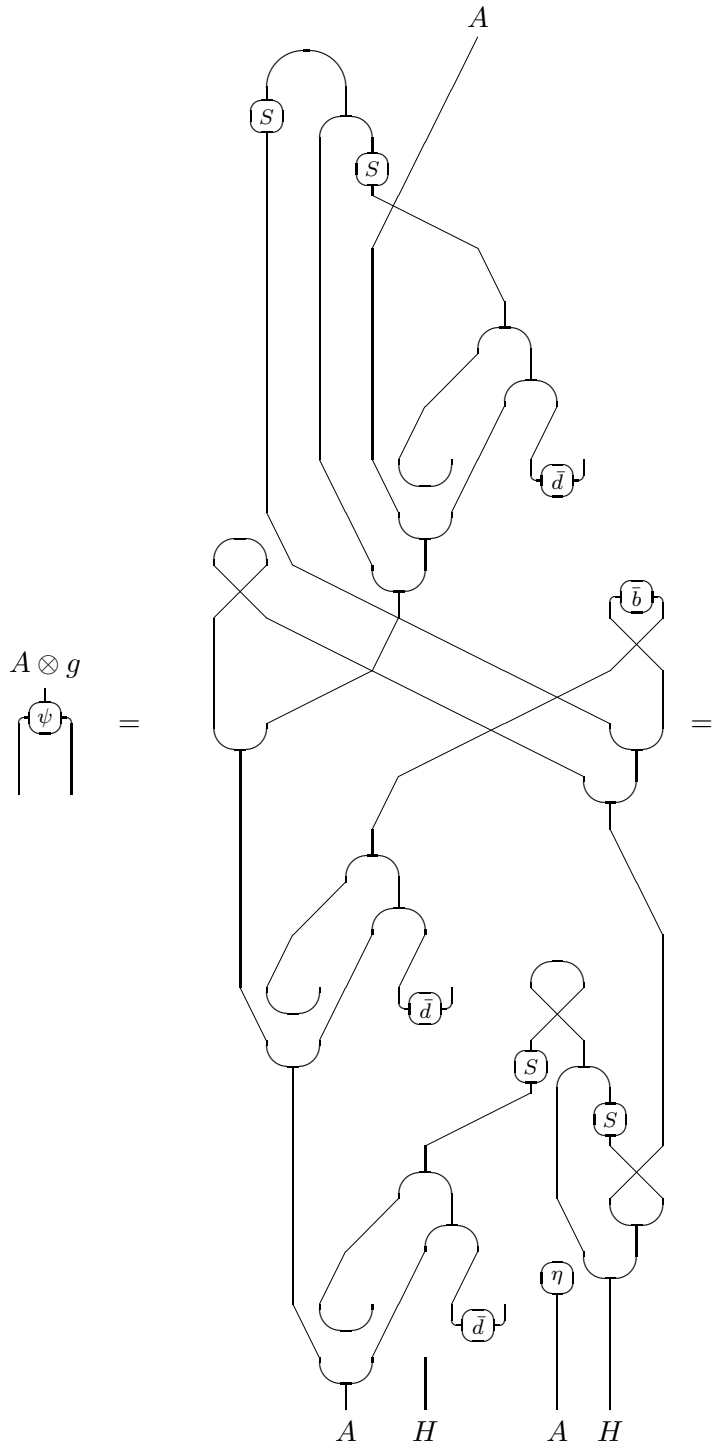
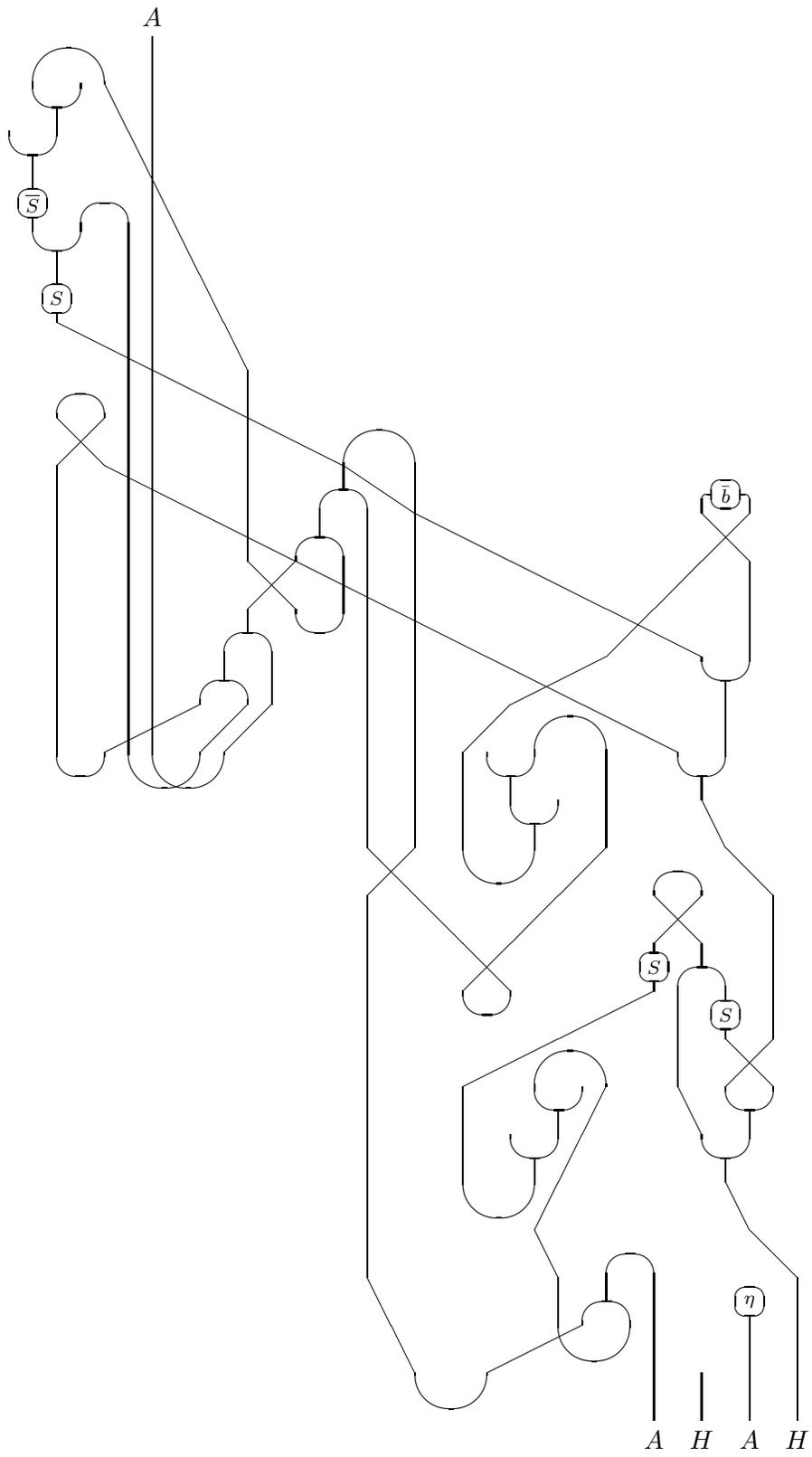
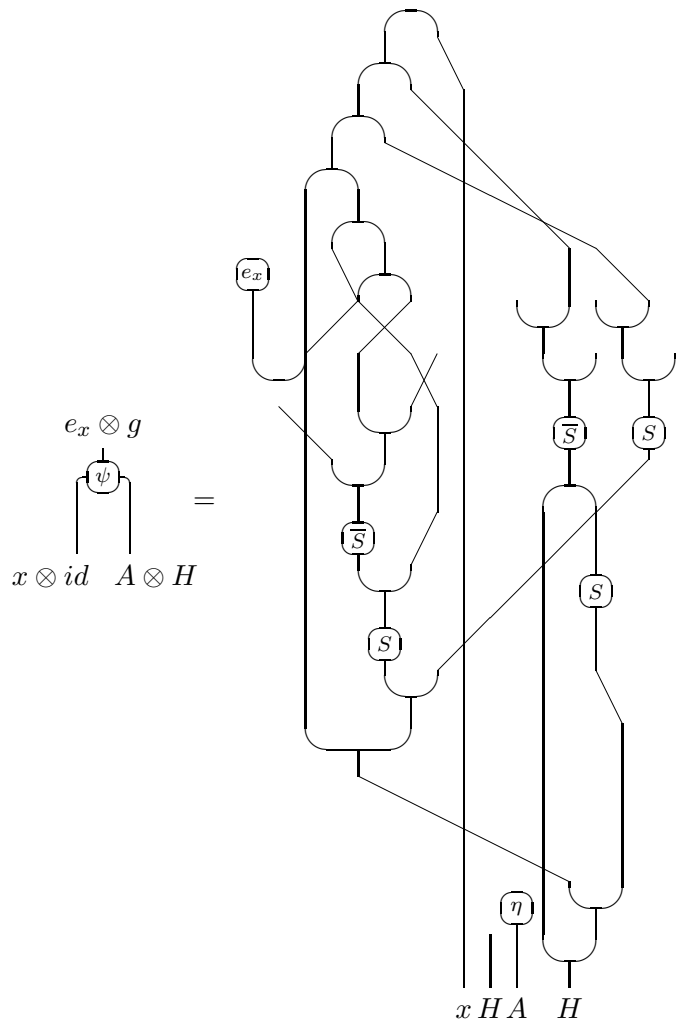
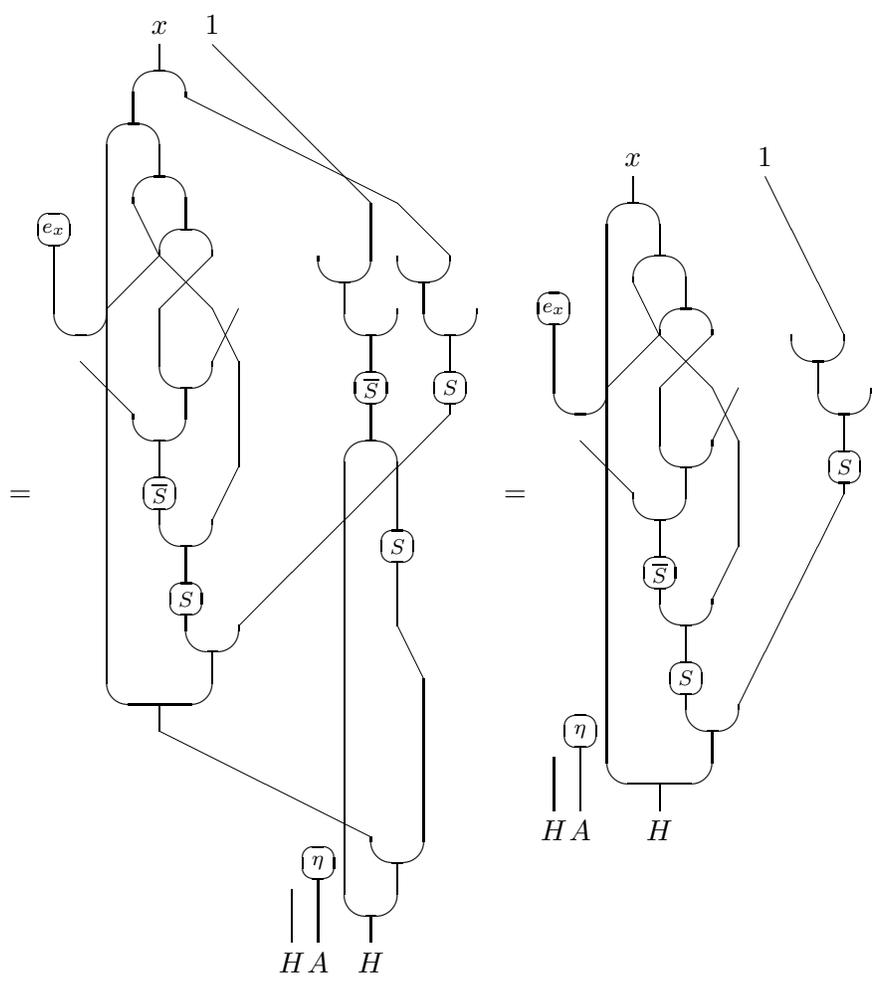
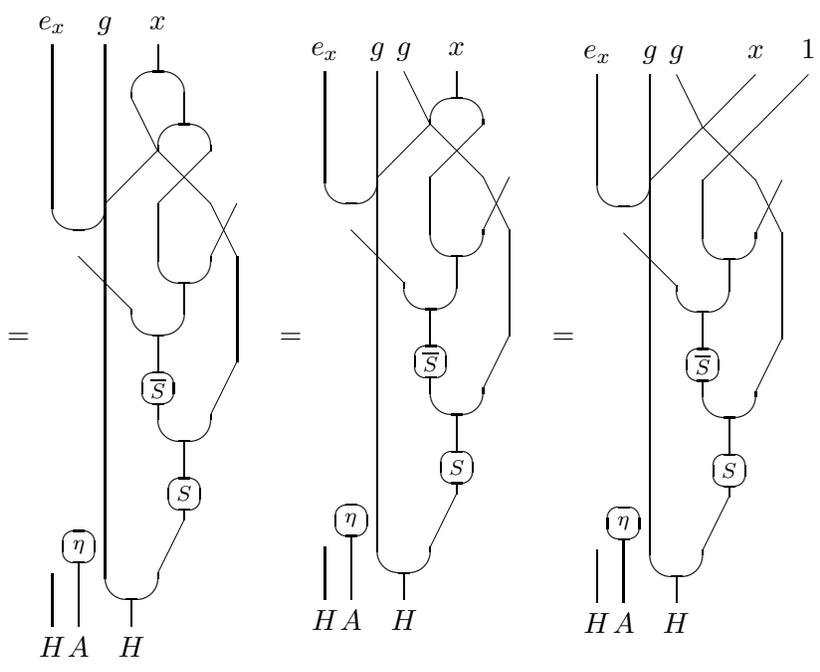
\[
\begin{tangle}
\step\object{A\otimes g}\\
\td \psi \\
\id\Step\id
\end{tangle}
\ \ \ \ = \ \ \ \
\begin{tangle}
\Step\Step\step\Step\Step\Step\object{A}\Step\object{}\\
\Step\step\Coev\Step\Step\Step\dd\step\\  
\Step\step\S\Step\cd\Step\dd\step\\  
\Step\step\id\Step\id\Step\S\step\dd\step\\    
\Step\step\id\Step\id\Step\ne1\step\step\se4\\   
\Step\step\id\Step\id\Step\id\Step\Step\d\\    

\Step\step\id\Step\id\Step\id\Step\Step\cd\\   

\Step\step\id\Step\id\Step\id\Step\step\ne2\step\cd\\   
 \Step\step\id\Step\id\Step\id\step\dd\Step\dd\step\dd\\   

\Step\step\id\Step\d\step\d\ev\dd\Step\coro {\bar d}\\    
\step\coev\d\Step\d\step\cu\\    
\step\XX\step\nw4\Step\cu\Step\step\Step\Step\ro {\bar b}\\    

\step\id\Step\nw4\Step\step\ne1\nw4\Step\step\Step\Step\XX\\   
\step\id\Step\Step\step\ne4\nw4\step\step\Step\nw4\Step\ne4\Step\id\\   
\step\cu\Step\Step\step\sw4\Step\nw4\Step\Step\cu\\    
\Step\id\Step\Step\step\dd\Step\Step\Step\step\cu\\   
\Step\id\Step\Step\cd\Step\Step\Step\Step\d\\   
\Step\id\Step\step\ne2\step\cd\Step\Step\Step\Step\d\\   
\Step\id\step\dd\Step\dd\step\dd\Step\Step\coev\Step\step\id\\   
\Step\d\ev\dd\Step\coro {\bar d}\Step\step\XX\Step\step\id\\    
\Step\step\cu\Step\Step\Step\Step\S\step\cd\Step\id\\    
\Step\Step\id\Step\Step\Step\Step\ne4\step\id\Step\S\Step\id\\   
\Step\Step\id\Step\Step\cd\Step\Step\id\Step\XX\\   
\Step\Step\id\Step\step\ne2\step\cd\Step\step\d\step\cu\\   
\Step\Step\id\step\dd\Step\dd\step\dd\Step\step\Q \eta\step\cu\\   
\Step\Step\d\ev\dd\Step\coro {\bar d}\Step\id\Step\id\\    
\Step\Step\step\cu\Step\id\Step\Step\step\id\Step\id\\    
\Step\Step\Step\object{A}\Step\step\object{H}\Step\Step\step\object{A}\Step\object{H}
\end{tangle}
\ \ = \]
\[
 \begin{tangle}
\Step\Step\Step\object{A}\\
\step\Coev\Step\Step\id\\   
\step\cu\step\d\step\id\\   
\cu\Step\step\d\id\\          
\step\O {\overline{S}}\Step\coev\step\id\d\\  
\step\cu\Step\id\step\id\step\d\\   
\Step\S\Step\step\id\step\id\Step\d\\  
\Step\nw4\Step\id\step\id\Step\step\d\\   
\Step\coev\step\id\step\id\nw4\Step\step\id\\   
\Step\XX\step\id\step\id\Step\Step\id\nw4\Step\step\Coev\\   
\Step\id\Step\nw4\id\step\id\Step\Step\id\Step\step\cd\step\se3\id
\Step\Step\Step\Step\Step\Step\ro {\bar b}\\ 
\Step\id\Step\step\id\step\id\Step\nw4\step\id\Step\cd\step\id\Step\id
\nw4\Step\Step\Step\Step\Step\step\XX\\        
\Step\id\Step\step\id\step\id\Step\Step\XX\nw4\step\id\step\id\Step\id
\Step\Step\nw4\Step\Step\Step\ne2\Step\id\\     
\Step\id\Step\step\id\step\id\Step\step\cd\step\ev\step\id\step\nw4\id
\Step\Step\Step\Step\step\ne2\step\se4\Step\id\\     
\Step\id\Step\step\id\step\id\Step\cd\step\id\Step\Step\id\Step\id
\Step\step\nw4\Step\step\ne4\Step\Step\cu\\       
\Step\id\Step\sw4\id\step\id\Step\step\ne2\ne2\Step\Step\id\Step\id
\Step\step\ne2\step\Coev\step\nw4\Step\Step\step\id\\     
\Step\ev\step\Ev\Ev\Step\Step\Step\Step\id\Step\id\Step\id
\step\cu\Step\step\id\Step\step\cu\\                  
\Step\Step\Step\Step\Step\Step\Step\step\id
\Step\id\Step\id\Step\cu\Step\id\Step\Step\nw1\\       
\Step\Step\Step\Step\Step\Step\Step\step\XX
\Step\Ev\Step\Step\ne2\Step\coev\Step\se2\\         
\Step\Step\Step\Step\Step\Step\Step\step\id

\Step\nw2\Step\Step\ne2\Step\Step\XX\step\Step\id\\      
\Step\Step\Step\Step\Step\Step\Step\step\id

\Step\Step\XX\Step\Step\Step\S\step\cd\Step\id\\   
\Step\Step\Step\Step\Step\Step\Step\step\id
\Step\Step\ev\Step\Step\step\ne4\step\id\Step\S\Step\id\\     
\Step\Step\Step\Step\Step\Step\Step\step\id

\step\step [3]\sw4\Step\Coev\step\Step\Step\id\Step\XX \\  
\Step\Step\Step\Step\Step\Step\Step\step\id

\Step\Step\id\step\Step\cu\dd\step\Step\d\step\cu\\   
\Step\Step\Step\Step\Step\Step\Step\step\id

\Step\Step\id\Step\cu\dd\step\Step\Step\cu\\          
 \Step\Step\Step\Step\Step\Step\Step\step\id
\Step\Step\Ev\Step\dd\step\Step\Step\Step\nw1\\       

\Step\Step\Step\Step\Step\Step\Step\step\id
\Step\Step\Step\step\nw1\step [2]\coev\Step\Step\nw2\\     

\Step\Step\Step\Step\Step\Step\Step\step\d
\step\Step\Step\Step\id \step \cd\step\id \Step\step\step\Q \eta\Step\id\\  

\Step\Step\Step\Step\Step\Step\Step\Step\d
\step\Step\sw4\step\step\Ev\Step\step\id\Step \Step\id\step\step\id \\ 

\Step\Step\Step\Step\Step\Step\Step\Step\step
\Ev\Step\Step\Step\Step\step\id\Step\id\Step\id \Step\id\\     
\Step\Step\Step\step\Step\Step\Step\Step\Step\Step\Step\Step\Step\Step
\object{A}\Step\object{H}\Step\object{A}\Step\object{H}\\
 \end{tangle}
\] and
\[
\begin{tangle}
\step[2]\object{e_{x}\otimes g}\\
\step \td \psi\\
\step \id\Step\id\\
 \object{x \otimes id}\Step\Step \object{A\otimes H}
\end{tangle}  \ \ \ = \ \ \
\begin{tangle}
\Step\Step\step \coev\\
\Step\Step\cd\step\d\\    
\Step\step\cd\step\nw2\step\id\\  
\Step\cd\step\nw4\Step\id\nw2\\  
\Step\id\step\cd\Step\step\id\step\nw4\nw2\\ 
\Q {e_x}\Step\id\step\d\cd\Step\id\Step\Step\id\step\nw2\\  
\id\Step\id\step\ne2\XX\Step\id\Step\cu\step\cu\\   
\ev\id\Step\id\Step\X\step\id\Step\step\cu\step\cu\\   
\step\nw2\id\Step\cu\step\id\step
\id\Step\step\step\O {\overline{S}}\Step\step\S\\       
\Step\id\step\cu\Step\id\step\id
 \Step\step\cd\step\ne2\\        
 \Step\id\Step\O {\overline{S}}
 \Step\dd\step\id
 \Step \step\id\step\ne2\id\\      

\Step\id\Step\cu\Step\id
\step\sw2\step\id\Step\S\\   
\Step\id\Step\step\S\Step\step\id
\ne2\step\step\id\Step\id\\       
\Step\id\Step\step\cu\step\id
\Step\step\id\Step\d\\  
\Step\Cu\Step\id
\Step\step\id\Step\step\id\\  
\Step\Step\nw4\Step\step\id
\Step\step\id\Step\step\id\\  
\Step\Step\Step\Step\id
\nw4\Step\id\Step\step\id\\     
\Step\Step\Step\Step\id
\Step\Q \eta \step\id\step\cu\\      
\Step\Step\Step\Step\id
\step\id\step\id\step\cu\\  
\Step\Step\Step\Step\object{x}\step
\object{H}\step\object{A}\Step\object{H}\\
\end{tangle}\]\[
 \;=\enspace
 \begin{tangle}
\step[4]\object x\step[2]\object 1\\
\Step\step\cd\step\nw2\step\\  
\Step\cd\step\nw4\Step\nw2\\  
\Step\id\step\cd\Step\step\step\nw4\nw2\\ 
\Q {e_x}\Step\id\step\d\cd\Step\Step\Step\id\step\nw2\\  
\id\Step\id\step\ne2\XX\Step\Step\cu\step\cu\\   
\ev\id\Step\id\Step\X\step\Step\step\cu\step\cu\\   
\step\nw2\id\Step\cu\step\id\step
\Step\step\step\O {\overline{S}}\Step\step\S\\       
\Step\id\step\cu\Step\id\step
 \Step\step\cd\step\ne2\\        
 \Step\id\Step\O {\overline{S}}\Step\dd\step
 \Step \step\id\step\ne2\id\\      
\Step\id\Step\cu\Step
\step\sw2\step\id\Step\S\\   
\Step\id\Step\step\S\Step\step
\ne2\step\step\id\Step\id\\       
\Step\id\Step\step\cu\step
\Step\step\id\Step\d\\  
\Step\Cu\Step
\Step\step\id\Step\step\id\\  
\Step\Step\nw4\Step\step
\Step\step\id\Step\step\id\\  
\Step\Step\Step\Step
\nw4\Step\id\Step\step\id\\     
\Step\Step\Step\Step
\Step\Q \eta \step\id\step\cu\\      
\Step\Step\Step\Step
\step\id\step\id\step\cu\\  
\Step\Step\Step\Step\object{}\step
\object{H}\step\object{A}\Step\object{H}\\
\end{tangle}\ \ = \ \
 \begin{tangle}
 \step[3]\object x\step[5]\object 1\\
\Step\cd\Step\step\step\d\\  
\Step\id\step\cd\Step\step\step\d\\ 
\Q {e_x}\Step\id\step\d\cd\Step\Step\d\\  
\id\Step\id\step\ne2\XX\Step\step\cu\\   
\ev\id\Step\id\Step\X\step\Step\cu\\   
\step\nw2\id\Step\cu\step\id\step
\Step\step\S\\       
\Step\id\step\cu\Step\id\step
 \Step\dd\\        
 \Step\id\Step\O {\overline{S}}\Step\dd
 \Step\dd\\      
\Step\id\Step\cu
\step\step\dd\\   
\Step\id\Step\step\S\Step\dd\\  
\step\Q \eta \step\id\Step\step\cu\step
\Step\step\\  
\id\step\id\step\Cu\Step\\  
\object{H}\step\object{A}\step\Step\object{H}\\
\end{tangle} \]\[
 \;=\enspace
 \begin{tangle}
\object{e_{x}}\Step\object{g}\Step\object{x}\\
\id\Step\id\step\cd\Step\\ 
\id\Step\id\step\d\cd\\  
\id\Step\id\step\ne2\XX\\   
\ev\id\Step\id\Step\X\\   
\step\nw2\id\Step\cu\step\id\\     
\Step\id\step\cu\Step\id\step\\    
 \Step\id\Step\O {\overline{S}}\Step\dd\\      
\Step\id\Step\cu\\   
\Step\id\Step\step\S\\  
\step\Q \eta \step\id\Step\dd\\  
\id\step\id\step\cu\\  
\object{H}\step\object{A}\Step\object{H}\\
\end{tangle}
 \;=\enspace
 \begin{tangle}
\object{e_{x}}\Step\object{g}\step\object{g}\Step\object{x}\\
\id\Step\id\step\d\cd\\  
\id\Step\id\step\ne2\XX\\   
\ev\id\Step\id\Step\X\\   
\step\nw2\id\Step\cu\step\id\\     
\Step\id\step\cu\Step\id\step\\    
 \Step\id\Step\O {\overline{S}}\Step\dd\\      
\Step\id\Step\cu\\   
\Step\id\Step\step\S\\  
\step\Q \eta \step\id\Step\dd\\  
\id\step\id\step\cu\\  
\object{H}\step\object{A}\Step\object{H}\\
\end{tangle}
 \;=\enspace
 \begin{tangle}
 \object{e_{x}}\Step\object{g}\step\object{g}\Step\step\object{x}\Step\object{1}\\
\id\Step\id\step\d\step\ne2\step\ne2\\  
\id\Step\id\step\ne2\XX\\   
\ev\id\Step\id\Step\X\\   
\step\nw2\id\Step\cu\step\id\\     
\Step\id\step\cu\Step\id\step\\    
 \Step\id\Step\O {\overline{S}}\Step\dd\\      
\Step\id\Step\cu\\   
\Step\id\Step\step\S\\  
\step\Q \eta \step\id\Step\dd\\  
\id\step\id\step\cu\\  
\object{H}\step\object{A}\Step\object{H}\\
\end{tangle}
\]
\\
\[
 \Step\Step\Step\;=\enspace
  \begin{tangle}
\Step\Q \eta \Step\Q \eta \\
 \id\Step\id\Step\id\\
\object{H}\Step\object{A}\Step\object{H}
  \end{tangle}
  \Step\; , \  but \enspace\Step\Step
   \begin{tangle}
  \object{e_{g}\otimes g}\Step\step\object{\eta_{D}}\\
 \id\Step\step\id\\
\object{x\otimes id}\Step\step\object{D}
  \end{tangle}
 \;= 0 \ .\enspace
\]
Thus $\psi$ is not trivial.
 \begin{picture}(5,5)
\put(0,0){\line(0,1){5}}\put(5,5){\line(0,-1){5}}
\put(0,0){\line(1,0){5}}\put(5,5){\line(-1,0){5}}
\end{picture}\\

{\bf Acknowledgement }  This work was supported by the National
Natural Science Foundation  (No.  19971074) and the fund of Hunan
education committee. Authors  thank the editors for valuable
suggestion and help.

\end{document}